\documentclass[12pt]{amsart}

\usepackage[margin = 1in]{geometry}

\usepackage{tikz}
\usepackage{tikz-cd}
\usepackage{url, hyperref}
\usepackage{float}

\tikzset{->-/.style={decoration={
  markings,
  mark=at position .45 with {\arrow{>}}},postaction={decorate}}}
  
\usetikzlibrary{shapes.arrows}
\usetikzlibrary{decorations.pathreplacing}
\usepackage[all]{xy}

\tikzset{->-/.style={decoration={
  markings,
  mark=at position .45 with {\arrow{>}}},postaction={decorate}}}

\usepackage{amsmath}
\usepackage{amssymb}
\usepackage{enumitem}
\usepackage{graphicx}
\usepackage{mathdots}
\usepackage{color}
\usepackage{diagbox}
\usepackage{array, makecell}
\usepackage{rotating}
\usepackage{amsmath}
\usepackage{tikz-cd}

\newcount\colveccount
\newcommand*\colvec[1]{
        \global\colveccount#1
        \begin{pmatrix}
        \colvecnext
}
\def\colvecnext#1{
        #1
        \global\advance\colveccount-1
        \ifnum\colveccount>0
                \\
                \expandafter\colvecnext
        \else
                \end{pmatrix}
        \fi
}

\usepackage{amsthm}

\theoremstyle{definition}
\newtheorem{definition}{Definition}[subsection]
\newtheorem{theorem}[definition]{Theorem}

\newtheorem{proposition}[definition]{Proposition}
\newtheorem{assumption}[definition]{Assumptions}

\newtheorem{corollary}[definition]{Corollary}
\newtheorem{lemma}[definition]{Lemma}

\newtheorem{remark}[definition]{Remark}

\newtheorem{question}[definition]{Question}

\def\bP{\mathbb{P}}
\def\bQ{\mathbb{Q}}

\def\bZ{\mathbb{Z}}

\def\cD{\mathcal{D}}
\def\cE{\mathcal{E}}

\def\cI{\mathcal{I}}

\def\cL{\mathcal{L}}
\def\cM{\mathcal{M}}

\def\cO{\mathcal{O}}

\def\barM{\overline{\mathcal{M}}}

\def\QEx{\mathrm{QEx}}

\DeclareMathOperator{\Gr}{Gr}

\DeclareMathOperator{\im}{im}

\DeclareMathOperator{\pr}{pr}

\newcommand{\wtE}{\widetilde{E}}
\newcommand{\wth}{\widetilde{h}}
\newcommand{\wtp}{\widetilde{p}}
\newcommand{\wtP}{\widetilde{P}}

\title{General curves on smooth quadric hypersurfaces}

 \author{Carl Lian}
\address{Washington University in St. Louis, Department of Mathematics, 1 Brookings Drive
\hfill \newline\texttt{}
 \indent  St. Louis, MO 63130} \email{{\tt clian@wustl.edu}}

\date{\today}

\usepackage{graphicx}

\begin{document}

\maketitle

\begin{abstract}
We prove a sharp existence theorem for non-degenerate maps from a general curve to a smooth quadric hypersurface, deforming in a smooth family of expected dimension. Along the way, we prove that if a general curve in projective space lies on a quadric (as dictated by the Maximal Rank Conjecture, now a theorem of E. Larson), then a general such quadric is smooth.
\end{abstract}

\setcounter{tocdepth}{1}

\section{Introduction}
\label{sec:intro}

\subsection{Existence}
\label{sec:intro_existence}
Throughout, we work over an algebraically closed field of characteristic 0.\footnote{The characteristic assumption will enter in two places: to invoke the main theorems of \cite{larson_reducible} in \S\ref{sec:inductive_smooth}, and to invoke generic smoothness in Proposition \ref{prop:existence_surjective_range}.} The study of curves and their deformations on higher-dimensional spaces is fundamental in geometry. In this paper, we consider the following question.

\begin{question}
\label{q:main}
Let $X$ be an algebraic variety, and let $C$ be a general curve of genus $g$. When does a map $f\colon C\to X$ of a given degree exist?
\end{question}

All curves in this paper are assumed to be projective and connected, with at worst nodes as singularities. The \emph{genus} of a curve refers to its arithmetic genus. By a \emph{general} curve, we mean a smooth curve lying in a non-empty open (but unspecified) subset of the moduli space of curves, $\cM_g$. 

When $X=\bP^r$ is a projective space, an answer to Question \ref{q:main} is given by \emph{Brill-Noether theory} \cite{gh}. Namely, a non-degenerate map $f\colon C\to\bP^r$ of degree $d$ exists if and only if
\[
\rho(g,r,d):=g-(r+1)(g-d+r)\geq 0.
\]
Recall that $f$ is \emph{non-degenerate} if its image is not contained in any hyperplane. It is natural to restrict to non-degenerate $f$, because degenerate $f$ can be viewed as maps to lower-dimensional projective spaces. 

The Gieseker-Petri theorem \cite{gieseker,eh_petri} furthermore asserts that the moduli space $\cM_d(C,\bP^r)^\circ$ of non-degenerate maps is smooth of dimension $\rho(g,r,d)+\dim\mathrm{Aut}(\bP^r)$, the expected. Equivalently, any non-degenerate map $f\colon C\to\bP^r$ out of a fixed, general curve has $H^1(C,f^*T_{\bP^r})=0$. Furthermore still, if $\rho(g,r,d)\geq 0$, then there exists a unique irreducible component $\barM_g(\bP^r,d)^{\mathrm{BN}}\subset\barM_g(\bP^r,d)$ dominating $\barM_g$, whose general member is a non-degenerate map out of a smooth curve \cite{eh_irreducible}. We refer to any point $f\colon C\to \bP^r$ on this component as a \emph{BN-curve} (even if $C$ is singular). 

As a landmark application, Eisenbud, Harris, and Mumford showed that the moduli space of curves $\barM_g$ is of general type when $g\geq 24$, by studying the \emph{Brill-Noether} divisor parametrizing curves admitting a non-degenerate map $f\colon C\to\bP^r$ of degree $d$ when $\rho(g,r,d)=-1$ \cite{hm,harris_odd,eh_kodaira}. This in particular requires the fact that no such map exists out of a general curve. In a recent breakthrough, Farkas-Jensen-Payne \cite{fjp2} showed the same for $\barM_{22}$ and $\barM_{23}$, using the divisors parametrizing curves admitting a non-degenerate map $f\colon C\to Q\subset\bP^6$ of degree $g+3$ to a quadric hypersurface. Once more, a crucial, and highly non-trivial piece of input is the non-existence of such a map out of a general curve, which was also proven independently by Liu-Osserman-Teixidor-Zhang \cite{lotz2}.

In this paper, we study the opposite question of \emph{existence} of maps $f\colon C\to Q\subset\bP^r$ to a smooth quadric hypersurface.

\begin{theorem}
\label{thm:main_existence}
Let $r\geq 4$ be an integer, and let $Q\subset\bP^r$ be a smooth quadric hypersurface. Let $C$ be a general curve of genus $g\ge0$. Then, a non-degenerate map $f:C\to Q$ of degree $d$ with $H^1(C,f^*T_Q)=0$ exists if and only if
\[
d\geq\max\left\{
\left\lceil\frac{r}{r+1}\cdot g+r\right\rceil,\,
g+\left\lfloor\frac r2\right\rfloor+1
\right\}.
\]
\end{theorem}

In particular, the moduli space of maps $C\to Q$ is smooth of expected dimension at $f$. As a consequence of Theorem \ref{thm:main_existence}, we obtain new examples of linear series on a general curve $C$ at which the Strong Maximal Rank Conjecture of Aprodu-Farkas \cite{af} holds, see \S\ref{sec:smr}.

See Remark \ref{rem:edge_cases} for a discussion of the edge cases $r=2,3$. A map $f:C\to Q$ is non-degenerate if it is so after post-composing with the inclusion $Q\subset\bP^r$. The condition
\[
d\geq \frac{r}{r+1}\cdot g+r
\quad\Longleftrightarrow\quad
\rho(g,r,d)\geq 0
\]
is necessary by the Brill-Noether theorem. If $H^1(C,f^*T_Q)=0$, then a map $f:C\to Q$ deforms in a smooth family of dimension
\[
\chi(C,f^*T_Q)=(r-1)(d-g+1),
\]
the expected. However, if $f$ is non-degenerate, then the $\frac{r(r+1)}{2}$-dimensional group of automorphisms of $Q$ acts freely at $f$, so we must have
\[
(r-1)(d-g+1)\geq \frac{r(r+1)}{2}
\quad\Longrightarrow\quad
d\geq g+\left\lfloor\frac{r}{2}\right\rfloor+1.
\]
Therefore, the ``if'' direction of Theorem \ref{thm:main_existence} is the interesting direction.

It remains an interesting question whether the assumption $d\geq g+\lfloor \frac r2\rfloor+1$ is necessary without the unobstructedness requirement $H^1(C,f^*T_Q)=0$. Teixidor i Bigas \cite{teixidor} showed that when $d\leq g+1$ and $r$ is arbitrary, no non-degenerate map $f$ exists. When $d=g+2$ and $r=4$, the same conclusion is established in  \cite[Theorem 2]{bafa} for all odd $g\ge11$. The results of \cite{fjp2,lotz2} prove non-existence when $(g,r,d)=(22,6,25)$ and $(23,6,26)$. While Theorem \ref{thm:main_existence} proves that the moduli space of non-degenerate maps $f\colon C\to Q$ has at least one component of expected dimension, we do not know whether it is the case that \emph{every} component of the space of non-degenerate maps has expected dimension.

We prove Theorem~\ref{thm:main_existence} in \S\ref{sec:existence}, by induction on $d$ and $g$. From maps out of curves of lower degree, we construct maps out of singular curves of the required degree, and then deform to maps out of smooth curves. This strategy has previously been developed, for example, in \cite{larson,cl_secant}. In this work, the inductive steps are relatively straightforward, but it is the base cases that are more difficult. We discuss this in the next section.

\begin{remark}\label{rem:edge_cases}
When $r=3$, the inequality on $d$ in Theorem \ref{thm:main_existence} is still necessary, but no longer sufficient. Indeed, when $g\ge 3$ is odd and $d=g+2$, then no non-degenerate map $f:C\to Q\cong\bP^1\times\bP^1\subset\bP^3$ can exist, because the gonality of a general curve of genus $g$ is $\left\lceil\frac{g+2}{2}\right\rceil$. However, it is a straightforward consequence of Brill-Noether theory for maps to $\bP^1$ that Theorem \ref{thm:main_existence} holds in all other cases when $r=3$.

Similarly, when $r=2$, Theorem \ref{thm:main_existence} is valid, as long as one restricts to even $d$. In this case, a multiple cover of a conic $Q\subset\bP^2$ of degree $\frac{d}{2}$ exists if and only if $d\ge g+2$.
\end{remark}

\subsection{Maximal rank}
\label{sec:intro_mr}

Rather than considering curves $f\colon C\to Q$ on quadrics, it is equally natural to consider quadrics containing a given curve. That is, let $f\colon C\to\bP^r$ be a map of degree $d$ out of a genus $g$ curve, and consider the restriction map
\[
\mu\colon H^0(\bP^r,\cO_{\bP^r}(2))\to H^0(C,\cO_C(2)).
\]
Then, the space of quadrics $Q \subset \bP^r$ containing the image of $f$ is identified with the projectivized kernel of $\mu$.

Suppose that $\rho(g,r,d) \geq 0$. The \emph{Maximal Rank Conjecture}, now a theorem of E. Larson \cite{larson_mr}, asserts that if $f \in \barM_g(\bP^r,d)^{\mathrm{BN}}$ is a general BN-curve, then $\mu$ is either injective or surjective\footnote{Larson proved the Maximal Rank Conjecture in its most general form, that, for a general BN-curve $f$, the restriction map $\mu$ is either injective or surjective for hypersurfaces in $\bP^r$ of \emph{arbitrary} degree. The case of quadrics had previously appeared in Ballico \cite{ballico}, and was proven independently by Jensen-Payne \cite{jp2}, and Liu-Osserman-Teixidor-Zhang \cite{lotz}. We do not attempt to survey the long history of the Maximal Rank Conjecture, which goes back at least to Severi \cite{severi}, here.}. We have $\dim H^0(\bP^r,\cO_{\bP^r}(2)) = \binom{r+2}{2}$, and a standard argument (see e.g. \cite[Remark 2.1]{jp2}) shows that $\dim H^0(C,\cO_C(2)) = 2d-g+1$. Therefore, if $f:C \to \bP^r$ is a general BN-curve and $\binom{r+2}{2} > 2d-g+1$, then the space of quadrics $Q \subset \bP^r$ on which $f$ lies is a projective space of dimension
\begin{equation}\label{eq:dim_ker}
\binom{r+2}{2}-1-(2d-g+1),
\end{equation}
the expected. We say that a map $f:C \to \bP^r$ \emph{satisfies maximal rank} if the map $\mu$ is either injective or surjective, depending on the sign of the quantity \eqref{eq:dim_ker}.

If $C$ is a general curve, $f:C \to \bP^r$ satisfies maximal rank and $Q$ is a \emph{smooth} quadric containing the image of $f$, then an incidence correspondence argument (see Proposition~\ref{prop:existence_surjective_range} for details) shows that $f:C \to Q$ is a map yielding the conclusion of Theorem \ref{thm:main_existence}. Indeed, in this case, the space of pairs of a map $f$ and a (varying) quadric $Q$ containing its image is smooth of expected dimension at $(f,Q)$. Projecting to the space of quadrics in $\bP^r$, by generic smoothness and the fact that any two smooth quadrics in $\bP^r$ are projectively equivalent, it follows that the space of maps to the \emph{fixed} quadric $Q$ is also smooth of expected dimension at $f$.

However, the smoothness of $Q$ is certainly not immediate. For example, if $C=Q\cap S\subset \bP^3$ is the complete intersection of a singular quadric and a smooth cubic, then $C$ is a canonically embedded curve of genus 4. Moreover, $Q$ is the unique quadric containing $C$, so the inclusion $C\subset \bP^3$ satisfies maximal rank, but lies on no smooth quadrics. Nevertheless, we show that if a \emph{general} BN-curve $f$ has image contained in a non-empty space of quadrics, as dictated by the Maximal Rank Conjecture, then a general such quadric is smooth. (The previous example shows that the open subset where this holds may be strictly smaller than the open subset where maximal rank holds.)

\begin{theorem}
\label{thm:main_smooth}
Let $g,r,d$ be integers with $\rho(g,r,d)\geq 0$ and $r\geq 3$. Let $f:C\to\bP^r$ be a general BN-curve of degree $d$ and genus $g$. 

Suppose that $\binom{r+2}{2}>2d-g+1$, so that the restriction map
\[
\mu\colon H^0(\bP^r,\cO_{\bP^r}(2))\to H^0(C,\cO_C(2))
\]
is surjective with non-zero kernel. Then, a general quadric hypersurface $Q\subset\bP^r$ containing the image of $f$ is smooth.
\end{theorem}

Kadik\"oyl\"u \cite[Theorem 1.1]{kadikoylu} proved that the space of quadrics of any given rank containing a general BN-curve $f$ has expected dimension if $g-d+r\le 1$, and Farkas-Rim\'anyi \cite[Theorem 1.9]{fr} obtained the same conclusion for quadrics of rank at most 6, when $(g,r,d)=(24,7,28)$. An equivalent formulation of Theorem~\ref{thm:main_smooth} is that the variety of quadrics of rank at most $r$ containing a general BN-curve $f$ has expected dimension, with no additional assumptions on $g,r,d$.

We prove Theorem~\ref{thm:main_smooth} in \S\ref{sec:smooth}, by following Larson's proof of the Maximal Rank Conjecture for quadrics \cite[\S4]{larson_mr}. Namely, we construct, by induction on $r$, singular curves in $\barM_g(\bP^r,d)^{\mathrm{BN}}$ where Theorem~\ref{thm:main_smooth} holds. In order to construct such curves, the main new step in our setting is to show that there exist elliptic normal curves on a smooth quadric $Q \subset \bP^r$ meeting $r+1$ given points on a smooth hyperplane section of $Q$ transversely, see Proposition~\ref{prop:existence_elliptic}. This is achieved by an enumerative calculation and subsequent transversality analysis, occupying most of \S\ref{sec:low_genus}.

In \S\ref{sec:inductive_existence}, we will see that Theorem~\ref{thm:main_smooth} provides sufficiently many base cases to deduce our main existence theorem, Theorem \ref{thm:main_existence}.

\subsection{Acknowledgments}

We gratefully acknowledge the support of NSF Grant DMS-2601731. We also thank Roya Beheshti, Alessio Cela, Gavril Farkas, David Jensen, Eric Jovinelly, Eric Larson, Montserrat Teixidor i Bigas, and Isabel Vogt for conversations related to this paper.

Due to physical constraints, the author had limited typing capacity during the preparation of this manuscript. He dictated portions of the paper, word by word, into an LLM via speech-to-text, with strict instructions to make no wording or mathematical changes.\footnote{The author is currently unaware of alternative speech-to-text software that is well-optimized for \LaTeX, but enthusiastically welcomes suggestions from the community.} An LLM was also used to produce TikZ code for the figure in \S\ref{sec:maps_hyperplane} from a hand-drawn diagram. 

After posting the first preprint version publicly, we performed an LLM audit on the manuscript for typographical and mathematical errors. Most resulting changes were minor, but two substantive mathematical issues were raised. First, the initial claim in Proposition \ref{prop:existence_chain} that the construction always produces a chain of lines was pointed out to be incorrect; we amended the statement and proof, with no changes to its application later in the paper. Second, a gap in \S\ref{sec:smr} was pointed out, which we filled by developing Theorem \ref{thm:existence_injective} and its proof in the appendix. Both fixes were implemented without further AI assistance.

All mathematical content and text is that of the author, with obvious intellectual debt to the work of many other mathematicians, including those named above.

\section{Curves of low genus}
\label{sec:low_genus}

In this section, we show that there exist rational and elliptic curves of degrees $r$ and $r+1$, respectively, interpolating through prescribed points on a hyperplane section of a smooth quadric hypersurface $Q \subset \bP^r$. These results will be used in \S\ref{sec:smooth} to construct singular curves satisfying the conclusion of Theorem \ref{thm:main_smooth}.

\subsection{Setup}
\label{sec:setup}

Throughout the section, fix:
\begin{itemize}
\item a smooth quadric hypersurface $Q\subset \bP^r$,
\item a hyperplane $H\subset \bP^r$ such that $Q':=Q\cap H$ is smooth, and
\item points $x_1,\ldots,x_{r+1}\in Q'$.
\end{itemize}

Consider the following generality assumptions on the points $x_i$.
\begin{assumption}
\label{assumption:generality}
\,
\begin{enumerate}
\item[(1)] $x_1,\ldots,x_{r+1}\in H$ are in linearly general position. That is, no subset of $s\le r$ of the $x_i$ lie on a linear subspace of dimension $s-2$.
\item[(2)] No two of the points $x_i$ lie on a line contained in $Q'$.
\item[(3)] For any collection of $r-1$ points among $x_1,\ldots,x_{r+1}$, the tangent hyperplanes $T_{x_i}Q'\subset H$ have no common point of intersection on $Q'$.
\end{enumerate}
\end{assumption}
(3) is equivalent to the statement that there does not exist a point $y\in Q'$ with the property that, for some $r-1$ points $x_i$ among $x_1,\ldots,x_{r+1}$, the points $y,x_i$ lie on a common line in $Q'$. It follows also that, for any collection of $r-2$ points among $x_1,\ldots,x_{r+1}$, there exist only finitely many (necessarily at most two) $y\in Q'$ lying on a line in $Q'$ with each of those $x_i$.

The main results of this section are as follows.

\begin{proposition}
\label{prop:existence_chain}
Suppose $r\geq 3$. Let $x_1,\ldots,x_r\in H$ be points satisfying Assumption \ref{assumption:generality}(1)-(2). Then, there exists a (0-pointed) stable map $h:R\to Q$ of degree $r$ with the following property. $R$ is a nodal curve of genus 0 containing distinct irreducible components $L_1,\ldots,L_r\subset R$, each mapping isomorphically to a distinct line $\ell_i\subset Q$ containing $x_i$, such that $\ell_i$ does not lie in $H$. Furthermore, $h$ is non-degenerate.
\end{proposition}

As we will see, one expects the construction we provide to produce a chain of lines $h:R=L_1\cup\cdots\cup L_r\to Q$. However, we will allow $R$ to acquire additional components, necessarily contracted by $h$.

\begin{proof}
We first construct the lines $\ell_i\subset Q$. Let $\ell_1\subset Q$ be any line through $x_1$, not contained in $H$; such a line exists because $Q\cap H$ is smooth, so $T_{x_1}Q\neq H$. Then, there exists a unique line $\ell_2\subset Q$ through $x_2$ meeting $\ell_1$, because $T_{x_2}Q$ does not contain $x_1$ by Assumption \ref{assumption:generality}(2), and thus does not contain $\ell_1$. Furthermore, because $T_{x_2}Q\cap \ell_1\neq x_1$, we have $\ell_2\not\subset H$. Continuing in this fashion, we find lines $\ell_1,\ldots,\ell_r\subset Q$ containing the points $x_1,\ldots,x_{r}$, respectively, such that $\ell_i\cap \ell_{i+1}$ is a point for each $i$. 

If the points $\ell_i\cap\ell_{i+1}$ are distinct, then we may define $h:R=L_1\cup\cdots\cup L_r\to Q$ to be a chain of lines mapping isomorphically to their images in $Q$. In general, we construct a stable map $h$ inductively as follows. First, let $h_1:R_1=L_1\to Q$ be the inclusion of $\ell_1$. Suppose now that $h_i:R_i\to Q$ is a stable map of degree $i$ mapping components $L_1,\ldots,L_i\subset R_i$ isomorphically to $\ell_1,\ldots,\ell_i$. If the point $\ell_i\cap\ell_{i+1}$ is identified with a smooth point $n_i\in L_i$, then attach $L_{i+1}$ to $R_i$ at $n_i$, and define $h_{i+1}$ by mapping $L_{i+1}$ isomorphically to $\ell_{i+1}$. If $n_i\in L_i$ is a node of $R_i$, then replace $n_i$ with a rational bridge $B_i$, and attach $L_{i+1}$ to a smooth point of $B_i$. Then, define $h_{i+1}$ by mapping $L_{i+1}$ isomorphically to $\ell_{i+1}$.

If Assumption \ref{assumption:generality}(1) holds, then the non-degeneracy is clear, because $x_1,\ldots,x_{r}$ span $H$, but no $\ell_i$ is contained in $H$.
\end{proof}

\begin{proposition}
\label{prop:existence_elliptic}
Suppose $r\geq 4$. Let $(E,p_1,\ldots,p_{r+1})\in \cM_{1,{r+1}}$ be a general pointed curve. Suppose that all of Assumptions \ref{assumption:generality}(1)-(3) hold. Then, there exists a map $h:E\to Q$ of degree $r+1$, such that:
\begin{itemize}
\item for all $i=1,\ldots,r+1$, we have $h(p_i)=x_i$,
\item the image of $h$ is not contained in $Q'=Q\cap H$.
\end{itemize}
\end{proposition}

Proposition \ref{prop:existence_elliptic} also holds when $r=3$, but our proof will use that $r\geq 4$, and we will only require the case $r\geq 4$.

\subsection{Elliptic curves on $Q$}
\label{sec:quasimaps}

The rest of this section is devoted to the proof of Proposition \ref{prop:existence_elliptic}. Because the points $x_i\in Q$ are not assumed to be general and we require a non-degeneracy condition on $h$, it does \emph{not} immediately suffice to show that the restricted tangent bundle of a general elliptic normal curve on $Q$ satisfies interpolation in the sense of \cite{larson,aly,lv}. We will instead detect the existence of a map $h$ via intersection theory on (a version of) the moduli space of quasimaps.\footnote{Our space will differ from the usual space of quasimaps \cite{cfkm} in that we will remove certain excess loci ``by hand,'' see Definition \ref{def:J}. For closely related virtual calculations on the actual space of quasimaps to a hypersurface, see \cite{ms}.} The setup is similar (but not identical) to that of \cite{lian_hyp}.

We first make the following basic observation: if the image of $h:E\to Q$ is not contained in $Q'=Q\cap H$, then the pullback of $H$ to $E$ must be precisely the divisor $p_1+\cdots+p_{r+1}$. In particular, we have 
\[
h^{*}\cO_{\bP^r}(1)=\cO(p_1+\cdots+p_{r+1})=:\cL.
\]
Therefore, we will consider the set of (quasi-)maps $h:E\to \bP^r$ with underlying line bundle $\cL$, and cut out loci corresponding to the conditions $h(p_i)=x_i$ and $\im(h)\subset Q$. 

\begin{definition}
\label{def:P}
Define
\[
\bP=\{h=[h_0:\cdots: h_r] \, | \, h_i\in H^0(E,\cL)\}=\bP(H^0(E,\cL)^{r+1})\cong\bP^{(r+1)^2-1}
\]
to be the space of non-zero $(r+1)$-tuples of sections in $H^0(E,\cL)$, up to simultaneous scaling.
\end{definition}
A general point of $\bP$ is a non-degenerate map $h:E\to \bP^r$, but $h$ may acquire base-points at the boundary of $\bP$.

We may choose coordinates $z_0,\ldots,z_r$ for $\bP^r$ in such a way that $H\subset \bP^r$ is cut out by $z_0$, and
\[
\begin{aligned}
x_1 &= [0:1:0:\cdots:0]\\
x_2 &= [0:0:1:\cdots:0]\\
&\vdots\\
x_r &= [0:0:0:\cdots:1],\\
x_{r+1} &= [0:1:1:\cdots:1].
\end{aligned}
\]
Write also $Q=V(F)\subset \bP^r$, where
\[
F=\sum_{1\leq i<j\leq r}a_{ij}z_iz_j+z_0\left(\sum_{0\leq k\leq r}b_kz_k\right).
\]
Note that $F$ has vanishing $z_i^2$-coefficient, for all $i=1,\ldots,r$, because $x_i\in Q$. Similarly, $x_{r+1}\in Q$ implies that
\(
\sum_{1\leq i<j\leq r}a_{ij}=0.
\)

\begin{definition}
\label{def:I}
Let $I\subset \bP$ be the subscheme cut out by the following linear equations:
\begin{itemize}
\item for all $i=1,\ldots,r$, we have $h_j(p_i)=0$ for all $j\in [0,r]\setminus i$,
\item we have $h_1(p_{r+1})=\cdots=h_{r}(p_{r+1})$.
\end{itemize}
\end{definition}

The first condition implies that $h_0$ vanishes at all of $p_1,\ldots,p_r$. Because $h_0\in H^0(E,\cO(p_1+\cdots+p_{r+1}))$, it follows that $h_0$ must automatically vanish at $p_{r+1}$ as well.

On the open locus of $\bP$ corresponding to maps $h:E\to \bP^r$, the subscheme $I\subset \bP$ restricts to the locus where $h(p_i)=x_i$ for all $i=1,\ldots,r+1$. Being the vanishing locus of $r(r+1)-1$ linear equations, $I$ is a linear space of dimension at least $r+1$. (In fact, it has dimension exactly $r+1$, but we will not need this.) Note furthermore that, for any $h\in I$, we have 
\begin{itemize}
\item $h_0\in H^0(E,\cL(-p_1-\cdots-p_{r+1}))\cong H^0(E,\cO)\subset H^0(E,\cL)$, and
\item for all $i=1,\ldots,r$, we have 
\[
h_i\in H^0(E,\cL(-p_1-\cdots-p_{i-1}-p_{i+1}-\cdots-p_{r}))\cong H^0(E,\cO(p_i+p_{r+1}))\subset H^0(E,\cL).
\]
\end{itemize}
Fix a choice of non-zero section $h_\bullet\in H^0(E,\cL(-p_1-\cdots-p_{r+1}))$.

We now wish to impose the requirement that maps $h:E\to \bP^r$ lie in $Q$, corresponding to the condition that
\[
F(h)=\sum_{1\leq i<j\leq r}a_{ij}h_ih_j+h_0\left(\sum_{0\leq k\leq r}b_kh_k\right)=0\in H^0(E,\cL^{\otimes 2}).
\]
Note, however, that $F(h)$ automatically vanishes along $p_1+\cdots+p_{r+1}$. Indeed, every individual summand vanishes at each of $p_1,\ldots,p_r$, and the vanishing at $p_{r+1}$ follows from $h_1(p_{r+1})=\cdots=h_{r}(p_{r+1})$ and \(
\sum_{1\leq i<j\leq r}a_{ij}=0\). Therefore, we obtain in fact an element
\[
F'(h):=\frac{F(h)}{h_\bullet}\in H^0(E,\cL)
\]
which may be regarded as a global section of the vector bundle $H^0(E,\cL)\otimes \cO_{I}(2)$ on $I$.

Instead of imposing the condition $F'(h)=0$ directly, we first impose the conditions that $F'(h)\in H^0(E,\cL)$ vanishes at each of $p_1,\ldots,p_r$ (but not $p_{r+1}$). Equivalently, $F(h)\in H^0(E,\cL^{\otimes 2})$ vanishes at each of $p_1,\ldots,p_r$ to order 2. For ease of notation, consider first the point $p_1$. Then, $F(h)$ vanishes to order 2 at $p_1$ if and only if
\[
h_1\left(b_1h_0+\sum_{1<j\leq r}a_{1j}h_j\right)\in H^0(E,\cL^{\otimes 2})
\]
does, because all other terms in $F(h)$ automatically vanish to order 2 at $p_1$. 

However, it is preferable to disallow $h_1(p_1)=0$, or else $h$ necessarily has a base-point at $p_1$. Note furthermore that the second factor $b_1h_0+\sum_{1<j\leq r}a_{1j}h_j$ automatically vanishes at $p_1$, so (after dividing by a choice of non-zero section of $\cO(p_1)$) may be viewed as an element of $H^0(E,\cL(-p_1))$. Globally, we may view $b_1h_0+\sum_{1<j\leq r}a_{1j}h_j$ as a section of the vector bundle $H^0(E,\cL(-p_1))\otimes \cO_{I}(1)$ on $I$. We therefore make the following definition.

\begin{definition}
\label{def:J}
Let $J\subset I$ be the subscheme cut out by the linear equations
\[
b_ih_0+\sum_{j\neq i}a_{ij}h_j=0\in H^0(E,\cL(-p_i)|_{p_i}),
\]
for all $i=1,\ldots,r$. By convention, we take $a_{ij}=a_{ji}$.
\end{definition}

$J$ is a linear space of expected dimension 1. Now, the pullback of $F'(h)$ to $J$ may be viewed as a global section $F''(h)$ of the line bundle $H^0(E,\cL(-p_1-\cdots-p_r))\otimes \cO_{J}(2)$. (Note at this stage that imposing the additional condition that $F'(h)$ vanishes at the remaining point $p_{r+1}$ does nothing.) 

\begin{definition}
\label{def:Z}
Let $Z\subset J$ be the subscheme cut out by the vanishing of $F''(h)$.
\end{definition}

Because $\cO_{J}(2)$ is ample, $Z$ is non-empty, regardless of whether it has the expected dimension of 0. (Its ``virtual degree'' is 2; see also Remark \ref{rem:vtev}.) If it is the case that every point $h\in Z$ is an $(r+1)$-tuple of sections of $\cL$ with no common vanishing point, and in addition $h_0\neq 0$, then $h$ defines a map $h:E\to Q$ satisfying the conclusion of Proposition \ref{prop:existence_elliptic}. It remains to verify these two conditions.

\begin{lemma}
\label{lem:bp}
Let $h\in Z$ be a point, and suppose that $h_i(p_i)=0$, for some $i=1,\ldots,r$. Let
\[
h'=[h'_0:\cdots: h'_r]
\]
be the tuple of sections of $\cL(-p_i)$ obtained by simultaneously twisting the sections $h_i$ down by $p_i$. Then, if the sections $h'_j$ do not all vanish at $p_i$, then the point $h'(p_i)\in Q$ lies in $T_{x_i}Q$. 
\end{lemma}

\begin{proof}
The $i$-th equation defining the subscheme $J\subset I$ is
\[
b_ih'_0(p_i)+\sum_{j\neq i}a_{ij}h'_j(p_i)=0,
\]
which is precisely the condition that $h'(p_i)$ lies in $T_{x_i}Q$.
\end{proof}

\begin{proposition}
\label{prop:bpf}
Let $h\in Z$ be a point. Then, the sections $h_j\in H^0(E,\cL)$ share no common vanishing point, except possibly if $h_0=0$ and the sections $h_1,\ldots,h_{r}$ vanish to order 1 at a unique point $p\neq p_1,\ldots,p_{r+1}$.
\end{proposition}

\begin{proof}
Let $b$ be the number of points among $p_1,\ldots,p_{r+1}$ on which the $h_j$ simultaneously vanish, and let $b'$ be the degree of the divisor away from $p_1,\ldots,p_{r+1}$, on which the $h_j$ simultaneously vanish. If the $h_j$ simultaneously vanish to order greater than 1 at any $p_i$, then we absorb the excess vanishing order into $b'$. Assume for sake of contradiction that $b+b'>0$, and let $h':E\to Q$ be the map of degree $d':=(r+1)-(b+b')$ obtained from $h$ by twisting down at base-points. 

First, consider the case $b=0$. We claim that $b'=1$ and $h_0=0$, so we are in the exceptional case allowed by the statement of the Proposition. Indeed, if $h_0$ vanishes at any point other than $p_1,\ldots,p_{r+1}$, then we must have $h_0=0$, because $h_0\in H^0(E,\cL(-p_1-\cdots-p_{r+1}))$. Now, $h'\colon E\to Q$ is a map of degree $d'\le r$, whose image contains $x_1,\ldots,x_{r+1}$. In particular, the linear span of the image of $h'$ is all of $H\cong\bP^{r-1}$, by Assumption \ref{assumption:generality}(1). Therefore (see e.g. Lemma \ref{lem:degree_elliptic}), we must have $d'=r$, hence $b'=1$.

Assume next that $b,d'>0$. Then, the image of $h'$ passes through $(r+1)-b\ge d'$ of the points $x_1,\ldots,x_{r+1}$, so in particular spans a linear space of dimension at least $d'-1$ by Assumption \ref{assumption:generality}(1). If $b'>0$, then these inequalities are strict, but the image of a map $h':E\to \bP^r$ of degree $d'$ spans a linear space of dimension at most $d'-1$. Therefore, we must have $b'=0$. Let $\Lambda\cong \bP^{d'-1}\subset \bP^r$ be the linear span of the image of $h'$.

Without loss of generality, let $p_1,\ldots,p_b$ be the points on which the $h_j$ simultaneously vanish. (While the point $p_{r+1}$ plays a distinguished role in the definition of $Z$, the argument that follows remains valid if $p_{r+1}$ is among the base-points.) Then, we have
\[
(h')^*\cO_{\bP^r}(1)=\cL(-p_1-\cdots-p_b)\cong\cO_E(p_{b+1}+\cdots+p_{r+1}).
\]
Furthermore, the map $h':E\to \Lambda$ is given by the complete linear series of the divisor $p_{b+1}+\cdots+p_{r+1}$. However, the images of the points $p_{b+1},\ldots,p_{r+1}$ do not lie on a hyperplane in $\Lambda$, which is a contradiction.

It remains to consider the case $d'=0$, so that $h':E\to Q$ is a constant map with image $y\in Q$. If $y\neq x_1,\ldots,x_{r+1}$, then the sections $h_0,\ldots,h_r$ must vanish simultaneously at $p_1,\ldots,p_{r+1}$ to order exactly 1, because $\deg(\cL)=r+1$. By Lemma \ref{lem:bp}, $y$ must furthermore lie in the intersection of the tangent hyperplanes to $Q$ at each of the points $x_1,\ldots,x_r$. We claim, however, that this is impossible. Indeed, if this were the case, then $x_1,\ldots,x_r$ would each lie on the tangent hyperplane $T_yQ$. However, the linear span of $x_1,\ldots,x_r$ is given by the hyperplane $H$, which cannot equal $T_yQ$ because $Q'=H\cap Q$ is smooth.

Next, suppose that $y$ is equal to one of $x_1,\ldots,x_r$; without loss of generality, $y=x_1$. Then, by Assumption \ref{assumption:generality}(2) and Lemma \ref{lem:bp}, $h_0,\ldots,h_r$ must vanish simultaneously at $p_2,\ldots,p_r$ to order at least 2, and at $p_{r+1}$ to order at least 1. However, we have
\[
2(r-1)+1>r+1=\deg(\cL),
\]
which is a contradiction. We obtain a contradiction similarly if instead $y=x_{r+1}$.
\end{proof}

\subsection{Maps to a hyperplane section}
\label{sec:maps_hyperplane}

In this section, we prove that any $h\in Z$ satisfies the conclusion of Proposition \ref{prop:existence_elliptic}. By Proposition \ref{prop:bpf}, it suffices to rule out two possibilities.
\begin{enumerate}
\item[(1)] $h$ defines a map $h\colon E\to Q'$ of degree $r+1$, with $h(p_i)=x_i$ and $h^*\cO_{\bP^r}(1)\cong\cL$.
\item[(2)] After twisting down at a base-point $p\neq p_i$, the point $h\in Z$ defines a map $h'\colon E\to Q'$ of degree $r$, satisfying $h'(p_i)=x_i$, with no condition on the underlying line bundle.
\end{enumerate}
The second case is easier, so is deferred to the end.

\begin{proposition}
\label{prop:map_hyperplane}
Let $(E,p_1,\ldots,p_{r+1})$ be a general curve of genus 1. There does not exist a map $h:E\to Q'$ of degree $r+1$, satisfying $h(p_i)=x_i$ for $i=1,\ldots,r+1$, and
\[
h^*\cO(1)\cong\cL=\cO_E(p_1+\cdots+p_{r+1}).
\]
\end{proposition}

We will prove Proposition \ref{prop:map_hyperplane} by degenerating $E$ to a singular curve. Namely, let $E_0=E\cup P$ be the union of $E$ with a smooth rational curve $P$, attached at a node $p_0\in E_0$. Let $p_1,\ldots,p_{r-2}\in E$ and $p_{r-1},p_r,p_{r+1}\in P$ be general smooth points. Then, fix a 1-parameter family of stable curves
\[
(E,p_1,\ldots,p_{r+1})\rightsquigarrow(E_0,p_1,\ldots,p_{r+1}).
\]
(See Remark \ref{rem:motivation_degeneration} for a discussion of motivation for this particular degeneration.)

If a map $h$ as in Proposition \ref{prop:map_hyperplane} exists, then by properness of the moduli space of stable maps, there exists a nodal curve $(\wtE_0,\wtp_1,\ldots,\wtp_{r+1})$ with stable contraction $(E_0,p_1,\ldots,p_{r+1})$, and a stable map $\wth:\wtE_0\to Q'$, such that $\wth(\wtp_i)=x_i$ for $i=1,\ldots,r+1$.

We may describe the curve $(\wtE_0,\wtp_1,\ldots,\wtp_{r+1})$ as follows. $\wtE_0$ contains two components (abusively denoted) $E,P$, mapping isomorphically to the corresponding components of $E_0$, which are connected by a tree $B$ of rational curves mapping to $p_0\in E_0$. Let $\wtp_0^+=E\cap B$ and $\wtp_0^-=B\cap P$. If $B$ is empty, then write $\wtp_0=\wtp_0^+=\wtp_0^-$ for the unique pre-image of $p_0$.

For each $i=1,\ldots,r-2$, the point $\wtp_i$ is either equal to $p_i\in E$, or lies on a tree of rational curves $T_i$ attached to $E$ at $p_i$. The same is true for $i=r-1,r,r+1$, with $E$ replaced by $P$. Let $\wtE$ be the union of $E$ and any trees of rational curves attached to $E$ at points other than $\wtp_0^+$. (This includes any trees $T_i$, where $i\leq r-2$.) Similarly, let $\wtP$ be the union of $P$ and any trees of rational curves attached to $P$ at points other than $\wtp_0^-$. Then, we have 
\[
\wtE_0=\wtE\cup_{\wtp_0^+}B\cup_{\wtp_0^-}\wtP.
\]

\begin{center}
\begin{tikzpicture}[
  x=1cm,y=1cm,
  line cap=round,line join=round,
  every node/.style={font=\small},
  point/.style={circle,fill,inner sep=1.2pt}
]

\coordinate (Nplus) at (-0.60,4.25);
\coordinate (Nminus) at (0.30,2.95);

\draw (0,7.00)
  .. controls (0,6.08333) and (-0.20,5.16667) .. (Nplus)
  node[pos=.15,point,label=above left:$p_1$] (Lone) {}
  node[pos=.32,point,label=right:$p_2$] {}
  node[pos=.78,point,label=right:$p_{r-2}$] {};

\draw (Nplus)
  .. controls (-0.64364,4.15) and (-0.68974,4.05) ..
  (-0.73813,3.95);

\draw (-0.82846,4.58) -- (Nplus) -- (Nminus) --
  (0.50769,2.65);

\draw (-1.20,3.50) -- (0.20,3.50);

\draw (0.42908,3.55)
  .. controls (0.38227,3.35) and (0.33924,3.15) ..
  (Nminus);

\draw (Nminus)
  .. controls (0.1234,2.05) and (0.0234,1.15) ..
  (0,0.25)
  node[pos=.30,point,label=right:$p_{r-1}$] {}
  node[pos=.65,point,label=right:$p_r$] {};

\draw (0,0.25) -- (0,-0.08);

\draw (Lone.center) -- ++(1.50,0)
  node[pos=.70,point,label=above:$\widetilde p_1$] {};

\draw (0,0.25) -- (1.50,0.25)
  node[pos=.70,point,label=above:$\widetilde p_{r+1}$] {};

\node[point,label=left:$\widetilde p_0^+$]
  at (Nplus) {};
\node[point,label=left:$\widetilde p_0^-$]
  at (Nminus) {};
\node[point,label=left:$p_{r+1}$]
  at (0,0.25) {};

\node at (0.20,5.55) {$\vdots$};
\node[above] at (0,7.00) {$E$};
\node[below] at (0,-0.08) {$P$};

\node[right] at (1.55,6.5875) {$T_1$};
\node[right] at (1.55,0.25) {$T_{r+1}$};

\node at (0.50,-0.90) {$\widetilde E_0$};

\draw (-2.48,6.95)--(-2.70,6.95)
      --(-2.70,4.23)--(-2.48,4.23);

\draw (-2.48,4.13)--(-2.70,4.13)
      --(-2.70,3.03)--(-2.48,3.03);

\draw (-2.48,2.93)--(-2.70,2.93)
      --(-2.70,0.12)--(-2.48,0.12);

\node[left] at (-2.80,5.59) {$\widetilde E$};
\node[left] at (-2.80,3.58) {$B$};
\node[left] at (-2.80,1.52) {$\widetilde P$};

\draw[->] (2.70,3.50) -- (4.90,3.50);

\coordinate (Rnode) at (6.40,3.50);

\draw (7.40,7.00)
  .. controls (7.40,5.83333) and (7.06667,4.66667) ..
  (Rnode)
  node[pos=.16,point,label=right:$p_1$] {}
  node[pos=.35,point,label=right:$p_2$] {}
  node[pos=.82,point,label=right:$p_{r-2}$] {};

\draw (Rnode)
  .. controls (6.31429,3.35) and (6.22306,3.20) ..
  (6.12633,3.05);

\draw (6.06864,4.00)
  .. controls (6.18698,3.83333) and (6.29744,3.66667) ..
  (Rnode);

\draw (Rnode)
  .. controls (7.06667,2.41667) and (7.40,1.33333) ..
  (7.40,0.25)
  node[pos=.25,point,label=right:$p_{r-1}$] {}
  node[pos=.50,point,label=right:$p_r$] {}
  node[pos=.75,point,label=right:$p_{r+1}$] {};

\node[point,label=left:$p_0$] at (Rnode) {};

\node at (7.40,4.98) {$\vdots$};
\node[above] at (7.40,7.00) {$E$};
\node[below] at (7.40,0.25) {$P$};
\node at (7.40,-0.90) {$E_0$};

\end{tikzpicture}
\end{center}

\begin{lemma}
\label{lem:degree_conic}
We have $\deg(\wth|_{\wtP})\geq 2$. If equality holds, then $\wtP=P$ maps isomorphically to the unique smooth conic on $Q'$ passing through $x_{r-1},x_r,x_{r+1}$.
\end{lemma}

\begin{proof}
Assume that $\deg(\wth|_{\wtP})\leq 2$, and consider the image $P_\bullet$ of $\wth|_{\wtP}$ in $Q'$. If $P_\bullet$ contains only lines, then by Assumption \ref{assumption:generality}(2), $P_\bullet$ must contain at least 3 lines, which is a contradiction. Therefore, $P_\bullet$ must be equal to a conic passing through all of $x_{r-1},x_r,x_{r+1}$. Such a conic spans a 2-plane $\Pi\subset H$, and by Assumption \ref{assumption:generality}(2), $\Pi$ intersects $Q'$ transversely at $P_\bullet$. Furthermore, $\wtP$ must map isomorphically to $P_\bullet$, as any degree 0 components will necessarily be unstable.
\end{proof}

\begin{lemma}
\label{lem:degree_elliptic}
Suppose that $\wth|_E$ is non-constant. Let $\Lambda\subset H\cong\bP^{r-1}$ be the linear span of the image of $\wth|_{\wtE}$. Then, we have
\[
\deg(\wth|_{\wtE})\geq\dim(\Lambda)+1.
\]
\end{lemma}

\begin{proof}
We proceed by induction on the number of components of $\wtE$. If $\wtE=E$ is smooth, then the claim follows by Riemann-Roch; a line bundle of degree $d$ has a $d$-dimensional space of sections, so $\dim(\Lambda)\leq d-1$. Now, attaching any rational component of degree $\delta$ to $\wtE$ increases $\dim(\Lambda)$ by at most $\delta$, so the desired inequality continues to hold.
\end{proof}

\begin{corollary}
\label{cor:degree_elliptic}
We have $\deg(\wth|_{\wtE})\geq r-2$.
\end{corollary}

\begin{proof}
If $\wth|_E$ is non-constant, then the claim follows immediately from Lemma \ref{lem:degree_elliptic}, because the points $x_1,\ldots,x_{r-2}$ span a linear space of dimension $r-3$. If instead $\wth|_E$ is constant, with image equal to $y\neq x_1,\ldots,x_{r-2}$, then for all $i=1,\ldots,r-2$, the tree $T_i$ contributes at least 1 to the degree of $\wth|_{\wtE}$, so we obtain the desired inequality once more. Finally, if, without loss of generality, $y=x_1$, then by Assumption \ref{assumption:generality}(2), for all $i=2,\ldots,r-2$, the tree $T_i$ contributes at least 2 to the degree of $\wth|_{\wtE}$. In particular, we have
\[
\deg(\wth|_{\wtE})\geq 2(r-3)\geq r-2,
\]
by the assumption that $r\geq 4$.
\end{proof}

Combining Lemma \ref{lem:degree_conic} and Corollary \ref{cor:degree_elliptic}, $\deg(\wth|_{\wtE})$ is equal to either $r-1$ or $r-2$.

\begin{proposition}
\label{prop:not_constant}
$\wth|_E$ cannot be constant.
\end{proposition}

\begin{proof}
Assume first that $\deg(\wth|_{\wtE})=r-1$, in which case we must have $\deg(\wth|_{\wtP})=2$ and $B=\emptyset$ by Lemma \ref{lem:degree_conic}. Write simply $\wtp_0:=\wtp_0^+=\wtp_0^-$. By Lemma \ref{lem:degree_conic}, $\wth$ maps $\wtP=P$ isomorphically to the unique conic $K\subset Q'$ through $x_{r-1},x_r,x_{r+1}$. Set $x_0=\wth(\wtp_0)\in K$.

If $\wth|_E$ is constant, then the image of $E$ must equal $x_0$, because $E\cap P=\wtp_0$. Because $(P,p_0,p_{r-1},p_r,p_{r+1})\in M_{0,4}$ is general, a degree $r-1$ map $\wth|_{\wtE}$ contracting $E$ to $x_0$ must exist for \emph{every} point of $K$, including $x_0=x_{r-1}$. By Assumption \ref{assumption:generality}(2), it must be the case that
\[
\deg(\wth|_{\wtE})\geq 2(r-2)>r-1,
\]
which is a contradiction.

Consider now the case $\deg(\wth|_{\wtE})=r-2$. Let $y\in Q'$ be the image of $E$ under $\wth$. Then, either $y$ is equal to one of $x_1,\ldots,x_{r-2}$, or is one of finitely many points lying on a line with each of $x_1,\ldots,x_{r-2}$, by Assumption \ref{assumption:generality}(3). In either case, the hypothesis on $y$ in Lemma \ref{lem:degree_3} below is satisfied, by Assumption \ref{assumption:generality}(2)-(3). On the other hand, $\wth|_{B\cup\wtP}$ is (possibly after contracting unstable components) a stable map of degree 3, sending $\wtp_0^+,\wtp_{r-1},\wtp_r,\wtp_{r+1}$ to $y,x_{r-1},x_r,x_{r+1}$, respectively. Lemma \ref{lem:degree_3} (with $\wtP$ in place of $\wtP\cup B$, and $\wtp_0$ in place of $\wtp_0^+$) shows that no such map can exist, so we have a contradiction once more.
\end{proof}

\begin{lemma}
\label{lem:degree_3}
Let $(P,p_0,p_{r-1},p_r,p_{r+1})\in M_{0,4}$ and $x_{r-1},x_r,x_{r+1}\in Q'$ be as above. (In particular, $(P,p_0,p_{r-1},p_r,p_{r+1})\in M_{0,4}$ is general.) Let $y\in Q'$ be an additional point, and assume that $y$ does not lie on a line in $Q'$ with any of $x_{r-1},x_r,x_{r+1}$.

Then, there do not exist a nodal curve $(\wtP,\wtp_0,\wtp_{r-1},\wtp_r,\wtp_{r+1})$ with stable contraction $(P,p_0,\allowbreak p_{r-1},p_r,p_{r+1})$, and a stable map $\wth:\wtP\to Q'$ of degree 3, such that $\wth(\wtp_i)=x_i$ for $i=r-1,r,r+1$, and $\wth(\wtp_0)=y$.
\end{lemma}

\begin{proof}
We may replace $P$ with a nodal rational curve $P_1\cup P_2$, where each component contains two of the marked points $p_i$. Let $\wth:\wtP_1\cup\wtP_2\to Q'$ be a stable map satisfying the required incidence conditions, such that $\wtP_i$ maps to $P_i$ under stable contraction. Then, one of $\wth|_{\wtP_1}$ and $\wth|_{\wtP_2}$ has degree at most $1$. However, this is impossible, because no two of the points $y,x_{r-1},x_r,x_{r+1}$ lie on a line in $Q'$.
\end{proof}

We assume henceforth that $\wth|_E$ is not constant.

\begin{lemma}
\label{lem:tree_degree}
For any $i=1,\ldots,r-2$, if $T_i$ is non-empty, then $\deg(T_i)=1$. Furthermore, there are no trees of rational curves attached to $E$ other than $T_i$ and $B\cup\wtP$.
\end{lemma}

\begin{proof}
Without loss of generality, assume for sake of contradiction that $\deg(T_1)\geq 2$. Suppose first that $\deg(\wth|_{\wtE})=r-2$. Then, delete $T_1,\wtp_1$ from $\wtE$. The new map $\wth|_{\wtE}$ has degree at most $r-4$, and passes through $r-3$ points of $Q'$ in linearly general position, so the span of its image has dimension at least $r-4$. Because $\wth|_E$ remains non-constant, this contradicts Lemma \ref{lem:degree_elliptic}.

Suppose instead that $\deg(\wth|_{\wtE})=r-1$, so $B=\emptyset$ and $\deg(\wth|_{\wtP})=2$. Write again $\wtp_0:=\wtp_0^+=\wtp_0^-$. By Lemma \ref{lem:degree_conic}, $\wtP=P$ maps isomorphically to the unique smooth conic $K\subset Q'$ passing through $x_{r-1},x_r,x_{r+1}$. Set $x_0=\wth(\wtp_0)\in K$. A general point of $K$ does not lie in the linear space spanned by $x_1,\ldots,x_{r-2}$, as $x_{r-1}$ does not, so by generality of $p_0\in P$, we may assume that $x_1,\ldots,x_{r-2},x_0$ span a linear space of dimension $r-2$. Then, $\wth|_{\wtE}$ is a map of degree $r-1$ whose image contains the $r-1$ points $x_1,\ldots,x_{r-2},x_0$. Applying an identical argument as in the case $\deg(\wth|_{\wtE})=r-2$ yields a similar contradiction to Lemma \ref{lem:degree_elliptic}.

Therefore, we conclude that any non-empty tree $T_i$, for $i=1,\ldots,r-2$, has $\deg(T_i)=1$. If there are any additional trees attached to $E$, then deleting those trees from $\wtE$ and repeating the arguments above similarly yields a contradiction.
\end{proof}

\begin{proposition}
\label{prop:line_bundle}
\,
\begin{enumerate}
\item[(a)] Suppose that $\deg(\wth|_{\wtE})=r-2$, and assume without loss of generality that the non-empty trees attached to $E$ are $T_1,\ldots,T_b$. Then, $b\leq r-4$, and $\deg(\wth|_E)=(r-2)-b$. Furthermore, we have
\[
\left.\wth^*\cO_{\bP^r}(1)\right|_E
\not\cong\cO_E(p_{b+1}+\cdots+p_{r-2}).
\]

\item[(b)] Suppose that $\deg(\wth|_{\wtE})=r-1$. As above, write $\wtp_0=\wtp_0^+=\wtp_0^-$ and $x_0=\wth(\wtp_0)$. Assume without loss of generality that the non-empty trees attached to $E$ are $T_1,\ldots,T_b$. Then, $b\leq r-3$, and $\deg(\wth|_E)=(r-1)-b$. Furthermore, we have
\[
\left.\wth^*\cO_{\bP^r}(1)\right|_E
\not\cong\cO_E(p_{b+1}+\cdots+p_{r-2}+\wtp_0).
\]
\end{enumerate}
\end{proposition}

\begin{proof}
Consider first (a). That $\deg(\wth|_E)=(r-2)-b$ follows from Lemma \ref{lem:tree_degree}. We have $b\leq r-4$ because $\wth|_E$ is assumed to be non-constant, and also cannot have degree 1. Now, the image of $\wth|_E$ contains $x_{b+1},\ldots,x_{r-2}$, so has linear span $\Lambda$ of dimension at least $r-3-b$. By Lemma \ref{lem:degree_elliptic}, equality holds. The final claim concerning $\wth^*\cO_{\bP^r}(1)$ follows from the same argument as in Proposition \ref{prop:bpf}; the map $\wth|_E$ cannot be given by the complete linear series of $p_{b+1}+\cdots+p_{r-2}$, because the images of the $p_i$ do not lie in a hyperplane.

The proof of (b) is identical, using (as in the proof of Lemma \ref{lem:tree_degree}) that $x_0$ does not lie in the linear span of $x_1,\ldots,x_{r-2}$.
\end{proof}

We are now ready to deduce the main results of the section.

\begin{proof}[Proof of Proposition \ref{prop:map_hyperplane}]
Let $\pi:\cE\to \Delta$ be a family of stable curves
\[
(E,p_1,\ldots,p_{r+1})\rightsquigarrow(E_0,p_1,\ldots,p_{r+1}),
\]
where $\Delta$ is the spectrum of a discrete valuation ring. Let $\cD\subset\cE$ be the divisor given by the sections $p_1,\ldots,p_{r+1}$. Assume for sake of contradiction that a map as in Proposition \ref{prop:map_hyperplane} exists. 

Then, by properness of the moduli space of stable maps, there exists (possibly after a finite base change) a birational modification $\beta:\widetilde{\cE}\to\cE$ and a family of stable maps $\wth:\widetilde{\cE}\to Q'$, sending the sections $p_i$ to the points $x_i\in Q'$. By abuse of notation, we denote the proper transform of $\cD\subset\cE$ also by $\cD\subset\widetilde{\cE}$. 

The restriction of $\wth$ to the special fiber of $\pi$ is a stable map $\wth:\wtE_0\to Q'$. With notation as above, we claim that
\[
\left.\wth^*\cO_{\bP^r}(1)\right|_E \cong
\begin{cases}
\cO_E(p_{b+1}+\cdots+p_{r-2}) &
\text{if } \deg(\wth|_{\wtE})=r-2,\\
\cO_E(p_{b+1}+\cdots+p_{r-2}+\wtp_0^+) &
\text{if } \deg(\wth|_{\wtE})=r-1,
\end{cases}
\]
which contradicts Proposition \ref{prop:line_bundle}.

The line bundle $\wth^*\cO_{\bP^r}(1)$ on $\wtE$ restricts to $\cO_{\wtE}(\cD)$ on the general fiber of $\cE$. Thus, the difference between these two line bundles is supported on the special fiber $\widetilde{\cE}_0$ of $\widetilde{\cE}$. More precisely, we have
\begin{equation}
\label{eq:line_bundles}
\wth^*\cO_{\bP^r}(1)\cong
\cO_{\wtE}(\cD)\otimes
\cO_{\wtE}\left(a_E\cdot E+\sum_{i=1}^{b}S_i+A\right),
\end{equation}
where $a_E\in\bZ$ is an integer, $A$ is a linear combination of components of $B\cup\wtP$, and $S_i$ is a linear combination of the components of the non-empty tree $T_i$, for $i=1,\ldots,b$. Furthermore, write
\[
S_i=a_iT_i'+T_i'',
\]
where $a_i\in\bZ$ and $T_i'\subset T_i$ is the unique component of $T_i$ meeting $E$. Here, the components of $\wtE_0$ are regarded as $\bQ$-Cartier (but not necessarily Cartier) divisors on $\wtE$.

We claim that $a_i=a_E$ for all $i=1,\ldots,b$. Indeed, restricting both sides of \eqref{eq:line_bundles} to $T_i$ and comparing degrees gives 
\[
1=1+(a_E-a_i)(E\cdot T_i').
\]
Here, we have used that any component of $T_i$ other than $T_i'$ has intersection number 0 against $\wtE_0$ and $\overline{\wtE_0\setminus T_i}$, so has intersection number 0 against $T_i$. Because $E\cdot T_i'>0$, the conclusion follows.

Now, restricting both sides of \eqref{eq:line_bundles} to $E$ gives
\[
\left.\wth^*\cO_{\bP^r}(1)\right|_E \cong
\cO_E(p_{b+1}+\cdots+p_{r-2}+a_0\cdot\wtp_0^+)
\]
for some integer $a_0$. Comparing degrees, we must have $a_0=0$ if $\deg(\wth|_{\wtE})=r-2$, and $a_0=1$ if $\deg(\wth|_{\wtE})=r-1$.
\end{proof}

\begin{proposition}\label{prop:map_hyperplane_bp}
Let $(E,p_1,\ldots,p_{r+1})$ be a general curve of genus 1. There does not exist a map $h\colon E\to Q'$ of degree $r$, satisfying $h(p_i)=x_i$ for $i=1,\ldots,r+1$.
\end{proposition}

\begin{proof}
We follow the proof of Proposition \ref{prop:map_hyperplane}. Replace $E$ with the same pointed curve $E_0=E\cup P$, and let
\[
\wtE_0=\wtE\cup_{\wtp_0^+}B\cup_{\wtp_0^-}\wtP
\]
and $\widetilde h\colon\wtE_0\to Q'$ be as before. Lemma \ref{lem:degree_conic} and Corollary \ref{cor:degree_elliptic} still hold, so we conclude in fact that $B=\emptyset$ and set $\wtp_0:=\wtp_0^+=\wtp_0^-$. Again, $\widetilde h$ maps $\wtP=P$ isomorphically to the unique conic $K\subset Q'$ through $x_{r-1},x_r,x_{r+1}$. Set $x_0=\widetilde h(\wtp_0)\in K$, which may be assumed not to lie in the span of $x_1,\ldots,x_{r-2}$.

A similar argument as in Proposition \ref{prop:not_constant} shows that $\widetilde h|_E$ cannot be constant. On the other hand, if $\widetilde h|_E$ is non-constant, then $\widetilde h|_{\wtE}$ has degree $r-2$, but its image has linear span of dimension at least $r-2$, contradicting Lemma \ref{lem:degree_elliptic}.
\end{proof}

\begin{proof}[Proof of Proposition \ref{prop:existence_elliptic}]
The non-emptiness of the subscheme $Z\subset\bP$ (Definition \ref{def:Z} and the subsequent discussion), combined with Proposition \ref{prop:bpf} and Proposition \ref{prop:map_hyperplane_bp}, shows that there exists a map $h:E\to Q$ of degree $r+1$ with $h(p_i)=x_i$ for all $i=1,\ldots,r+1$. By construction of $Z$, we have $h^*\cO_{\bP^r}(1)\cong\cL$, even if the image of $h$ is contained in $Q'$. Then, Proposition \ref{prop:map_hyperplane} shows that no such map $h:E\to Q'$ exists, so any $h\in Z$ satisfies the conclusion of Proposition \ref{prop:existence_elliptic}.
\end{proof}

\begin{remark}\label{rem:vtev}
The locus on $I$ (Definition \ref{def:I}) where $h_0=0$ is a hyperplane, and is in particular ample. Therefore, the fact that any point $h\in Z$ has $h_0\neq 0$ (Proposition \ref{prop:map_hyperplane}) shows that $\dim(Z)=0$. Furthermore, $Z$ is the complete intersection in $\bP$ of hyperplanes and a single quadric, so has length 2. (Neither statement is logically necessary for us.)

This is consistent with Buch-Pandharipande's calculation of ``virtual Tevelev degrees'' of quadric hypersurfaces \cite[Theorem 1.4]{bp}.\footnote{On the other hand, our calculation takes place on a different moduli space, so neither calculation immediately implies the other.} Indeed, the virtual (Gromov-Witten) number of maps $h:E\to Q$ satisfying $h(p_i)=x_i$ is
\[
\mathrm{vTev}^{Q}_{g=1,d=r+1,n=r+1}
=
\begin{cases}
r-2, & r\text{ even},\\
r+1, & r\text{ odd}.
\end{cases}
\]
(Note that \cite[Theorem 1.4]{bp} uses the letter $r$ to denote the dimension of $Q$, whereas here $r$ is the dimension of the ambient projective space.) On the other hand, the virtual number of maps $h:E\to Q'$ satisfying $h(p_i)=x_i$ is
\[
\mathrm{vTev}^{Q'}_{g=1,d=r+1,n=r+1}
=
\begin{cases}
r-4, & r\text{ even},\\
r-1, & r\text{ odd}.
\end{cases}
\]
In particular, in both cases, we have
\begin{equation}\label{eq:vtev}
\mathrm{vTev}^{Q}_{1,r+1,r+1}
-
\mathrm{vTev}^{Q'}_{1,r+1,r+1}
=2.
\end{equation}
Therefore, one expects the number of maps $h$ as in Proposition \ref{prop:existence_elliptic} to be 2.

However, because both counts are virtual, it is not clear how to deduce either Proposition \ref{prop:existence_elliptic} or Proposition \ref{prop:map_hyperplane} directly from \eqref{eq:vtev}. The points $x_i$ not being assumed to be general, the enumerativity arguments of e.g. \cite{lian_pand,bllrst} do not immediately apply.
\end{remark}

\begin{remark}\label{rem:motivation_degeneration}
We comment on the motivation for choosing the degenerate curve $E_0=E\cup P$. In the setting of Proposition \ref{prop:map_hyperplane}, naive parameter counts lead one to expect that the only maps $h\colon E_0\to Q'$ satisfying the incidence conditions $h(p_i)=x_i$ have $\deg(h|_{E})=r-1$ and $\deg(h|_P)=2$. Furthermore, if, as above, $x_0$ is the image of $p_0$ under the unique isomorphism sending $P$ to the conic $K\subset Q'$ containing $x_{r-1},x_r,x_{r+1}$, then we require $h(p_0)=x_0$. This is consistent with the equality of virtual counts
\[
\mathrm{vTev}^{Q'}_{1,r+1,r+1}=\mathrm{vTev}^{Q'}_{1,r-1,r-1},
\]
which holds by \cite[Theorem 1.4]{bp}. 

The advantage of passing from $h$ to $h|_{E}$ is that the degree decreases by 2, but the dimension of the span of the image decreases by only 1. Thus, the argument constraining the underlying line bundle in the proof of Proposition \ref{prop:map_hyperplane} (and earlier in the proof of Proposition \ref{prop:bpf}) becomes available. Note that we did not need to prove that the only maps $h\colon E_0\to Q'$ are those described above, but constrained the possible stable maps sufficiently to apply this argument. The only potential obstruction was the possibility that a stable map out of a model of $E_0$ is constant on the main component $E$, but this is ruled out by placing sufficiently many (namely, $r-2$) marked points on $E$.
\end{remark}

\section{Smooth quadrics containing curves}
\label{sec:smooth}

In this section, we prove Theorem \ref{thm:main_smooth} by induction on $r$. 

\subsection{Base cases}
\label{sec:base_cases}

\begin{proposition}
\label{prop:smooth_base_cases}
Theorem \ref{thm:main_smooth} holds when $r=3$.
\end{proposition}

\begin{proof}
The required inequalities $d\geq \frac{3}{4}g+3$ and $10>2d-g+1$ hold for only finitely many pairs $(d,g)$, namely:
\[
(d,g)\in\{(3,0),(4,0),(4,1),(5,2),(6,4)\}.
\]
We will consider each of these cases separately. (All of these cases are covered by \cite[Theorem 1.1]{kadikoylu}, but we give a self-contained proof.) Write $z_0,z_1,z_2,z_3$ for coordinates on $\bP^3$. Because the locus where Theorem \ref{thm:main_smooth} holds is open, it suffices to produce \emph{one} example of a map $f$ satisfying maximal rank, and lying on \emph{one} smooth quadric.

\underline{$(d,g)=(3,0)$:} Any non-degenerate $f:\bP^1\to\bP^3$ is linearly equivalent to the map
\[
f([s:t])=[s^3:s^2t:st^2:t^3].
\]
The space of quadric forms pulling back to 0 under $f$ is spanned by $z_0z_3-z_1z_2$, $z_0z_2-z_1^2$, $z_1z_3-z_2^2$, so $f$ satisfies maximal rank. The quadric cut out by $z_0z_3-z_1z_2$ is smooth.

\underline{$(d,g)=(4,0)$:} The non-degenerate map $f:\bP^1\to\bP^3$ given by
\[
f([s:t])=[s^4:s^3t:st^3:t^4]
\]
satisfies maximal rank, as the only quadric pulling back to 0 is $z_0z_3-z_1z_2$. This quadric is smooth.

\underline{$(d,g)=(4,1)$:} For any genus 1 curve $E$, any non-degenerate map $f:E\to\bP^3$ is defined by the complete linear series underlying a line bundle $\cL$ of degree 4 on $E$. For any such $\cL$, write $\cL=\cL_1\otimes\cL_2$, where $\cL_1\not\cong\cL_2$ are line bundles of degree 2; this is possible because the squaring map $\operatorname{Pic}^2(E)\to\operatorname{Pic}^4(E)$ is finite. By the base-point-free pencil trick \cite[Chapter III, \S 3]{acgh}, the map
\[
\otimes:H^0(E,\cL_1)\otimes H^0(E,\cL_2)
\longrightarrow H^0(E,\cL)
\]
is injective, hence an isomorphism. In particular, any $f:E\to\bP^3$ lies on a smooth quadric, namely that given by the Segre embedding underlying $\cL=\cL_1\otimes\cL_2$.

\underline{$(d,g)=(5,2)$:} For any genus 2 curve $C$, any non-degenerate map $f:C\to\bP^3$ is defined by the complete linear series underlying a line bundle $\cL$ of degree 5 on $C$. We may write $\cL=K_C\otimes\cM$, where $K_C$ is the canonical line bundle on $C$ and $\cM$ has degree 3. If $f$ is general, then $\cM$ is base-point-free. Indeed, if $\cM$ were not base-point-free, then $\cM\cong K_C(p)$ for some $p\in C$. However, the space of such line bundles has dimension 1 in the 2-dimensional Picard variety $\mathrm{Pic}^3(C)$. Arguing as in the previous case, the base-point-free pencil trick shows that the map
\[
\otimes:H^0(C,K_C)\otimes H^0(C,\cM)
\longrightarrow H^0(C,\cL)
\]
is an isomorphism, so $f$ lies on a smooth quadric.

\underline{$(d,g)=(6,4)$:} Let $C$ be a general curve of genus 4. Then, any non-degenerate map $f:C\to\bP^3$ of degree 6 must be a canonical embedding. Furthermore, $C$ has exactly two \emph{distinct} (up to projective linear transformation) maps $f_i:C\to\bP^1$ of degree 3 \cite[Main Theorem IIb]{gh}. Write $\cL_i=f_i^*\mathcal{O}_{\bP^1}(1)$. By the base-point-free pencil trick, the map
\[
\otimes:H^0(C,\cL_1)\otimes H^0(C,\cL_2)
\longrightarrow H^0(C,\cL_1\otimes\cL_2)
\]
is injective. Because $\cL_1\otimes\cL_2$ has degree 6, we also have $h^0(C,\cL_1\otimes\cL_2)\leq 4$, so in fact, $\cL_1\otimes\cL_2\cong f^*\mathcal{O}_{\bP^3}(1)$, and $f$ lies on a smooth quadric.
\end{proof}

\subsection{Inductive steps}
\label{sec:inductive_smooth}

Assume throughout this section that $r\ge 4$.

\begin{lemma}
\label{lem:degenerate_rational}
\,
\begin{enumerate}
\item[(a)] Let $\delta$ be an integer with $2 \leq \delta \leq r-2$, and suppose that Theorem \ref{thm:main_smooth} holds in $\bP^\delta$ for a rational normal curve $f \colon \bP^1 \hookrightarrow \bP^\delta$. Then, Theorem \ref{thm:main_smooth} holds in $\bP^{r-1}$ for the composition $f \colon \bP^1 \hookrightarrow \bP^\delta \subset \bP^{r-1}$, in the sense that the composition
\[
H^0(\bP^{r-1},\cO_{\bP^{r-1}}(2))
\to H^0(\bP^\delta,\cO_{\bP^{\delta}}(2))
\to H^0(\bP^1,\cO_{\bP^{\delta}}(2))=H^0(\bP^1,\cO_{\bP^1}(2\delta))
\]
is surjective, and a general element of the kernel cuts out a smooth quadric.
\item[(b)] The conclusion of (a) also holds for $\delta=0,1$, in the following sense. If $\bP^1 \hookrightarrow \bP^{r-1}$ and $x \hookrightarrow \bP^{r-1}$ are the inclusions of a line and a point, respectively, then the maps $H^0(\bP^{r-1},\cO_{\bP^{r-1}}(2))\to H^0(\bP^1,\cO_{\bP^{1}}(2))$ and $H^0(\bP^{r-1},\cO_{\bP^{r-1}}(2)) \to \cO(2)|_x$ are surjective, and a general element of either kernel cuts out a smooth quadric.
\end{enumerate}
\end{lemma}

\begin{proof}
(a) follows from the fact that $H^0(\bP^{r-1},\cO(2)) \to H^0(\bP^\delta,\cO(2))$ is surjective, and any smooth quadric in $\bP^\delta$ lifts to a smooth quadric in $\bP^{r-1}$. (Strictly speaking, Theorem \ref{thm:main_smooth} is only stated when $\delta \geq 3$, but when $\delta=2$, the hypothesis that $H^0(\bP^2,\cO_{\bP^{2}}(2)) \to H^0(\bP^1,\cO_{\bP^{2}}(2))$ is surjective with a smooth conic in the kernel clearly holds.) (b) is straightforward.
\end{proof}

\begin{lemma}
\label{lem:point_general}
Let $C''$ be a smooth curve, and let $f'':C''\to H\cong\bP^{r-1}$ be any non-degenerate map, where $r\geq 4$. Let $Q'\subset H$ be a smooth quadric containing the image of $C''$. Then, a general collection of points $x_1,\ldots,x_{r+1}\in C''$ satisfies all properties of Assumption \ref{assumption:generality}.
\end{lemma}

Here, despite the fact that $f''$ is not assumed to be injective, we harmlessly identify points $x_i\in C''$ with their images in $Q'$.

\begin{proof}
The space of ordered tuples of points on $C''$ is irreducible, and all of the conditions in Assumption \ref{assumption:generality} are open, so it suffices to show that there exists \emph{one} collection of points $x_i$, for which each of the assumptions (1)-(3) holds separately. For (1), one may choose the points $x_i\in C''$ one at a time, in such a way that $x_i$ does not lie in the linear span of any subset of the previously chosen points of size at most $r-1$. This is possible because no proper linear subspace of $H$ contains the image of $C''$. Similarly, for (2), one may choose the points $x_i\in C''$ one at a time, in such a way that $x_i$ does not lie in the tangent hyperplane to $Q'$ at any of the previously chosen points.

Finally, to arrange for (3), we may assume that (1) already holds. Then, the condition that the tangent hyperplanes at $x_1,\ldots,x_{r-1}$ do not intersect on $Q'$ is equivalent to the condition that the hyperplane $\Lambda \subset H$ spanned by $x_1,\ldots,x_{r-1}$ intersects $Q'$ transversely. Indeed, viewing $Q' \subset H$ as the space of isotropic vectors for a non-degenerate, symmetric bilinear form, we have
\[
y \in x_1^\perp,\ldots,x_{r-1}^\perp \Longleftrightarrow \langle x_1,\ldots,x_{r-1}\rangle\subset y^\perp 
\]
which is impossible if $y\in Q'$ and $\Lambda \cap Q'$ is smooth.

A general hyperplane $\Lambda\subset H$ has smooth intersection with $Q'$. Furthermore, a general hyperplane $\Lambda\subset H$ has the property that its intersection with $C''$ contains some collection of $r-1$ points in linearly general position; for example, choose $x_1,\ldots,x_{r-1}\in C''$ in linearly general position, and take $\Lambda$ to be their linear span. Therefore, a general hyperplane $\Lambda\subset H$ satisfies both properties simultaneously. In particular, there exist $x_1,\ldots,x_{r-1}\in C''$ in linearly general position, with the property that their linear span $\Lambda$ intersects $Q'$ smoothly. This implies more generally that, for general points $x_1,\ldots,x_{r+1} \in C''$, the linear span of any $r-1$ of the $x_i$ intersects $Q'$ smoothly. Indeed, the locus of such $x_i$ is given by the intersection of finitely many non-empty open subsets of the space of ordered tuples of points on $C''$.
\end{proof}

The remainder of the proof closely follows \cite[\S4]{larson_mr}. It suffices to show that there exists one BN-curve $f:C\to\bP^r$ satisfying maximal rank, and lying in at least one smooth quadric. We produce such a curve by induction on $r$, by attaching a curve $h:C'\to Q\subset\bP^r$ to a curve $f'':C''\to H\cong\bP^{r-1}\subset\bP^r$, where $Q$ is a smooth lift of a quadric in $H$ containing the image of $f''$, and $h$ is a map constructed in Proposition \ref{prop:existence_chain} or \ref{prop:existence_elliptic}. 

In an effort to make the exposition as self-contained as possible, we repeat many aspects of the proof from \cite[\S4]{larson_mr}. However, we emphasize that the key new input, from the previous section, is that $h$ may be taken to lie on a \emph{smooth} quadric hypersurface containing the image of $C''$.

\begin{proof}[Proof of Theorem \ref{thm:main_smooth}]
Fix a hyperplane $H\cong\bP^{r-1}\subset\bP^r$. We assume that the conclusion of Theorem \ref{thm:main_smooth} holds in $\bP^{r-1}$.

Suppose first that \underline{$d\geq g+r$ and $g\geq r-1$}. Write $d''=d-r$ and $g''=g-(r-1)\geq 0$. Note that $d''\geq g''+(r-1)$; in particular, $\rho(g'',r-1,d'')\geq 0$. Let $C''$ be a general curve of genus $g''$, and let $f'':C''\to H$ be a general BN-curve of degree $d''$. Note that
\begin{equation}\label{eq:kernel_equal}
\binom{r+1}{2}-(2d''-g''+1)
=
\binom{r+2}{2}-(2d-g+1),
\end{equation}
so by the inductive hypothesis, we may assume that Theorem \ref{thm:main_smooth} holds for $f''$. Fix a smooth quadric $Q'\subset H$ containing the image of $f''$ and a smooth quadric $Q\subset\bP^r$ with the property that $Q\cap H=Q'$. By Lemma \ref{lem:point_general}, we may choose general points $x_1,\ldots,x_r\in C''$ satisfying Assumptions \ref{assumption:generality}(1)-(2). 

By Proposition \ref{prop:existence_chain}, there exists a map $h:C'=R\to Q$ from a union of smooth rational curves, mapping components $L_i\subset R$ isomorphically onto a line $\ell_i\subset Q$ meeting $H$ transversely at $x_i$. Let $C=C'\cup C''$ be the curve obtained by gluing $C'$ and $C''$ along $x_1,\ldots,x_r$ at nodes, and let $f:C\to Q\subset\bP^r$ be the map obtained by gluing $h$ and $f''$. Then, $f$ is a stable map of degree $d$ and genus $g$.

The $r$ points $x_1,\ldots,x_r$ are in linearly general position on $H$, and can therefore be regarded as general points on $\bP^r$. Furthermore, the space of stable maps $\barM_0(\bP^r,r)$ is irreducible \cite{kp}, so $h:C'\to\bP^r$ is a BN-curve. By \cite[Theorem 1.6]{larson_reducible}, it follows that $f\in\barM_g(\bP^r,d)$ is a BN-curve. We claim furthermore that $f$ satisfies maximal rank, cf. \cite[Proposition 2.2]{larson_mr}. 

Indeed, we have a commutative square
\[
\begin{tikzcd}
H^0(\cI_f(2)) \arrow[r, "j"] \arrow[d, hook] &
H^0(\cI_{f''}(2)) \arrow[d, hook] \\
H^0(\cO_{\bP^r}(2)) \arrow[r] &
H^0(\cO_H(2))
\end{tikzcd}
\]
defined by restricting quadrics from $\bP^r$ to $H$, where $H^0(\cI_f(2))$ and $H^0(\cI_{f''}(2))$ are the subspaces of quadrics in their respective projective spaces vanishing on $C$ and $C''$. Then, $j$ is injective, because a non-zero element of the kernel would necessarily cut out the union of $H$ with another hyperplane, which cannot vanish on $C'$ by the non-degeneracy statement of Proposition \ref{prop:existence_chain}. 

On the other hand, we have
\[
\begin{aligned}
\dim H^0(\cI_f(2))
&\geq \binom{r+2}{2}-(2d-g+1)\\
&= \binom{r+1}{2}-(2d''-g''+1)\\
&= \dim H^0(\cI_{f''}(2)),
\end{aligned}
\]
because $f''$ satisfies maximal rank. It follows that $H^0(\cI_f(2))$ also has expected dimension. Because the image of $f$ is contained in a smooth quadric by construction, $f$ satisfies the conclusion of Theorem \ref{thm:main_smooth}.

Next, suppose that \underline{$d\geq g+r$ and $g\leq r-2$}. We also exclude the case $(g,d)=(1,r+1)$, and defer this case to the end. Write again $d''=d-r$. Let $C''$ be a smooth curve of genus 0, and let $f'':C''\to H$ be a general map of degree $d''$, which is either a general BN-curve (if $d''\geq r-1$) or a rational normal curve in a proper linear subspace $\Lambda\subset H$ (if $d''\leq r-2$). If $d''\geq r-1$, then write $\Lambda=H$. We have
\[
\binom{r+1}{2}-(2d''-0+1)
=
\binom{r+2}{2}-(2d-g+1)+(r-(g+1)).
\]

By the inductive hypothesis and Lemma \ref{lem:degenerate_rational}, we may assume that Theorem \ref{thm:main_smooth} holds for $f''$, and choose a smooth quadric $Q'\subset H$ containing the image of $f''$. If $2\leq d''\leq r-2$, then we may also arrange for $Q'$ to intersect the linear span $\Lambda$ of $C''$ in a smooth quadric (by specifying first a smooth quadric $Q'\cap \Lambda$ containing the image of $f''$, and then lifting to a smooth quadric $Q'\subset H$). 

We claim that we may choose general points $x_1,\ldots,x_{g+1}\in C''$ and $x_{g+2},\ldots,x_{r}\in Q'$, satisfying the following conditions.
\begin{enumerate}
    \item[(i)] $x_1,\ldots,x_r$ satisfy Assumptions \ref{assumption:generality}(1)-(2). 
    \item[(ii)] $x_{g+2},\ldots,x_{r}$ impose independent conditions on quadrics in $H$ containing $C''$. That is, the space $H^0\bigl(\cI_{f''\cup\{x_{g+2},\ldots,x_{r}\}}(2)\bigr)$ of quadrics in $H$ vanishing on $C''$ and $x_{g+2},\ldots,x_{r}$ has dimension
\[
\binom{r+1}{2}-(2d''+1)-(r-g-1),
\]
the expected.
\end{enumerate}

Indeed, arguing as in the proof of Lemma \ref{lem:point_general}, we may begin by choosing general points $x_1,\ldots,x_{g+1}\in C''$, satisfying Assumptions \ref{assumption:generality}(1)-(2). Here, if $d''\leq r-2$, then we initially view $C''$ as lying on the smooth quadric $Q'\cap\Lambda\subset\Lambda\cong\bP^{d''}$. It is at this step that we must exclude the case $(g,d)=(1,r+1)$. In this case, the image of $C''$ is a line on which we must choose two points $x_1,x_2$, so Assumption \ref{assumption:generality}(2) cannot hold. In all other cases, either $d''\geq 2$, in which case $Q'\cap\Lambda\subset\Lambda$ is a smooth quadric in a linear space of dimension at least 2, or $g=0$, in which case we only choose $x_1\in C''$. 

If $d''\leq r-2$ and $x_1,\ldots,x_{g+1}$ are in linearly general position in $\Lambda\cong\bP^{d''}$, then they are also in linearly general position in $H\cong\bP^{r-1}$, because
\[
g+1\leq d''+1=\dim(\Lambda)+1.
\]
Similarly, if $x_1,\ldots,x_{g+1}$ satisfy Assumption \ref{assumption:generality}(2) in $\Lambda$, then the same is true in $H$. Now, the remaining general points $x_{g+2},\ldots,x_r\in Q'$ may be chosen, one at a time, in such a way that Assumptions \ref{assumption:generality}(1)-(2) still hold. Furthermore, condition (ii) can be guaranteed, by choosing $x_i$ not to lie in the base locus of $H^0(\cI_{f''\cup\{x_{g+2},\ldots,x_{i-1}\}}(2))$.

Choose now a smooth quadric $Q\subset\bP^r$ with $Q\cap H=Q'$. Finally, define a map $h:C'=R\to Q$ as in the previous case, along with a map $f:C=C'\cup C''\to Q$. Once more, $f$ is a BN-curve\footnote{If $d=r$ (and necessarily $g=0$), then $f|_{C''}$ is constant, so one needs to contract $C''$ in order to obtain a stable map.}, and we claim that $f$ satisfies maximal rank. Indeed, consider the restriction map
\[
j:H^0(\cI_f(2))
\longrightarrow
H^0\bigl(\cI_{f''\cup\{x_{g+2},\ldots,x_{r}\}}(2)\bigr).
\]
By condition (ii) above, the dimension of the source of $j$ is at least that of the target. On the other hand, a non-zero element of the kernel of $j$ is a quadric given by the union of $H$ and a hyperplane vanishing on the image of $h$, which is impossible. We conclude that $f$ satisfies the conclusion of Theorem \ref{thm:main_smooth}.

Finally, suppose that \underline{$d<g+r$}. Write $d''=d-r-1$ and $g''=g-r-1$. Then, exactly as in \cite[\S4]{larson_mr}, we have
\[
\begin{aligned}
g'' = \rho(g,r,d) + (r+1)(g+r-d-1) &\geq 0,\\
\rho(g'',r-1,d'') = \rho(g,r,d) + (g+r-d-1) &\geq 0.
\end{aligned}
\]
Let $C''$ be a general curve of genus $g''$, and let $f'':C''\to H$ be a general BN-curve of degree $d''$. The relation \eqref{eq:kernel_equal} holds, so Theorem \ref{thm:main_smooth} holds for $f''$.

Choose smooth quadrics $Q'\subset H$ containing the image of $f''$, and $Q\subset\bP^r$ such that $Q\cap H=Q'$. By Lemma \ref{lem:point_general}, we may choose $x_1,\ldots,x_{r+1}\in C''$ satisfying all properties of Assumption \ref{assumption:generality}. Let $h:C'=E\to Q$ be a map as in Proposition \ref{prop:existence_elliptic}. Define $f:C=C'\cup C''\to Q\subset\bP^r$ by gluing along the nodes $x_1,\ldots,x_{r+1}$. The image of $h$ meets $H$ transversely at the distinct points $x_i$, which may be regarded as general points in $H\cong\bP^{r-1}$ by linear generality. By \cite[Theorem 1.9]{larson_reducible}, $f$ is a BN-curve.\footnote{Here, we also require the inequality $d''+(r+1)\ge g''+r$, which follows from $\rho(g,r,d)\ge0$ and $\binom{r+2}{2}>2d-g+1$.}  By a similar argument as in the previous cases, $f$ satisfies maximal rank, and therefore the conclusion of Theorem \ref{thm:main_smooth}.

The remaining case $(g,d)=(1,r+1)$ can be dealt with in a similar way to the case $d<g+r$. Namely, for general points $x_1,\ldots,x_{r+1}\in H$, there exists a smooth quadric $Q'\subset H$ containing each of the $x_i$, and such that the $x_i$ impose independent conditions on quadrics in $H$. Choose a smooth lift $Q\subset\bP^r$ with $Q\cap H=Q'$, and a map $h\colon E\to Q$ as in Proposition \ref{prop:existence_elliptic}. Then, $H^0(\mathcal{I}_h(2))$ is identified with the space of quadrics in $H$ containing $x_1,\ldots,x_{r+1}$, which has dimension
\[
\binom{r+1}{2}-(r+1)=\binom{r+2}{2}-2(r+1).
\]
\end{proof}

\section{Curves on smooth quadrics}
\label{sec:existence}

\subsection{Setup}
\label{sec:existence_setup}

In this section, we prove Theorem~\ref{thm:main_existence}. Let $r\geq 4$ be an integer, and fix a smooth quadric hypersurface $Q\subset\bP^r$. For $g,r,d$ satisfying the hypothesis of Theorem~\ref{thm:main_existence}, write $\QEx(g,r,d)$ for the statement that there exists a non-degenerate map $f:C\to Q$ out of a general, smooth curve $C$ of genus $g$ with $H^1(C,f^*T_Q)=0$. We say that $f$ is ``nice'' if it is non-degenerate and $H^1(C,f^*T_Q)=0$.

The notion of niceness makes sense for any stable map $f:C\to Q$, and by deformation theory, if there exists a nice stable map out of a single nodal curve, then the same is true out of a general, smooth curve. We will prove Theorem~\ref{thm:main_existence} by induction on $d$ and $g$, by constructing nice maps out of singular curves from nice maps out of smooth curves of lower degree. Base cases will be afforded by Proposition~\ref{prop:existence_surjective_range}.

\subsection{Base cases}
\label{sec:existence_base_cases}

\begin{proposition}
\label{prop:existence_surjective_range}
Suppose that $\rho(g,r,d)\geq 0$ and $\binom{r+2}{2}>2d-g+1$, so that Theorem~\ref{thm:main_smooth} holds. Then, $\QEx(g,r,d)$ holds.
\end{proposition}

The inequalities $\rho(g,r,d)\geq 0$ and $\binom{r+2}{2}>2d-g+1$ imply that $d\ge g+\lfloor \frac r2\rfloor+1$, consistent with the necessary conditions for $\QEx(g,r,d)$ to hold.

\begin{proof}
Fix a general curve $C$ of genus $g$, and denote by $\cM_d(C,Q)\subset\cM_d(C,\bP^r)$ the spaces of degree $d$ maps $f:C\to Q$ and $f:C\to\bP^r$, respectively. Then, the open subset of $\cM_d(C,\bP^r)$ parameterizing non-degenerate maps is non-empty and smooth of dimension $\chi(C,f^*T_{\bP^r})$. Let $\cM\subset\cM_d(C,\bP^r)$ be the non-empty open subset consisting of non-degenerate maps satisfying maximal rank.

Consider the incidence correspondence
\[
\begin{tikzcd}
& Y = \{(f,F) \mid \operatorname{im}(f) \subset V(F)\}
\arrow[dl, "\pr_1"']
\arrow[dr, "\pr_2"] & \\
\cM && \bP H^0(\cO_{\bP^r}(2))
\end{tikzcd}
\]
of pairs of a map $f:C\to\bP^r$ along with a quadric containing the image of $f$. Then, $\pr_1$ is a projective bundle whose fibers have dimension $\left(\binom{r+2}{2}-1\right)-(2d-g+1)$, so $Y$ is non-empty and smooth of dimension
\[
\chi(C,f^*T_{\bP^r})+\left(\binom{r+2}{2}-1\right)-(2d-g+1).
\]

By Theorem~\ref{thm:main_smooth}, the image of $Y$ under $\pr_2$ meets the locus of smooth quadrics in $\bP H^0(\cO_{\bP^r}(2))$. In fact, the image of $Y$ meets the entire locus of smooth quadrics, because $\mathrm{Aut}(\bP^r)$ acts transitively on smooth quadrics. The locus of smooth quadrics in $\bP H^0(\cO_{\bP^r}(2))$ is dense, so by generic smoothness, the fiber of $\pr_2$ over any smooth quadric is smooth of dimension
\[
\chi(C,f^*T_{\bP^r})-(2d-g+1).
\]
In particular, the fiber $\pr_2^{-1}([Q])$ is non-empty and smooth of this dimension. On the other hand, $\pr_2^{-1}([Q])$ is identified with an open subset of $\cM_d(C,Q)$, consisting of non-degenerate maps $f:C\to Q$.

The expected dimension of $\cM_d(C,Q)$ is
\[
\begin{aligned}
\chi(C,f^*T_{Q})&=\chi(C,f^*T_{\bP^r})-\chi(C,\cO_{\bP^r}(2)|_C)\\
&=\chi(C,f^*T_{\bP^r})-(2d-g+1)\\
&=\dim\left(\pr_2^{-1}([Q])\right),
\end{aligned}
\]
so any $f\in\pr_2^{-1}([Q])$ defines a nice $f:C\to Q$. In particular, $\QEx(g,r,d)$ holds.
\end{proof}

\subsection{Inductive steps}
\label{sec:inductive_existence}

The following two lemmas are standard. For a systematic study of restricted tangent bundles of rational curves on hypersurfaces, see \cite{mioranci}.

\begin{lemma}
\label{lem:line_tangent}
Let $h:\bP^1\hookrightarrow Q$ be the inclusion of a line. Then, 
\[
h^*T_Q \cong \cO_{\bP^1}\oplus \cO_{\bP^1}(1)^{\oplus(r-3)}\oplus \cO_{\bP^1}(2).
\]
\end{lemma}

\begin{proof}
Let $\bP^1\hookrightarrow \bP^1\times\bP^1\hookrightarrow Q$ be the sequence of inclusions given by a ruling line, followed by a 2-dimensional quadric obtained by intersecting $Q$ with $r-3$ general hyperplanes containing that line. Then, the short exact sequence 
\[
0\to N_{\bP^1/(\bP^1\times\bP^1)}\to N_{\bP^1/Q}\to N_{(\bP^1\times\bP^1)/Q}|_{\bP^1}\to 0
\]
yields a splitting $N_{\bP^1/Q}\cong \cO_{\bP^1}\oplus \cO_{\bP^1}(1)^{\oplus(r-3)}$, because $\operatorname{Ext}^1(\cO_{\bP^1}(1),\cO_{\bP^1})=0$. The short exact sequence \[
0\to T_{\bP^1}\cong\cO_{\bP^1}(2)\to h^*T_Q\to N_{\bP^1/Q}\to 0
\]
also splits, yielding the desired conclusion.
\end{proof}

\begin{lemma}
\label{lem:conic_tangent}
Let $h:\bP^1\hookrightarrow Q$ be the inclusion of a smooth conic, given by the transverse intersection of $Q$ with a 2-plane. Then, 
\[
h^*T_Q \cong \cO_{\bP^1}(2)^{\oplus(r-1)}.
\]
\end{lemma}

\begin{proof}
The smooth conic in $Q$ is cut out by $r-2$ general hyperplanes in $\bP^r$, so $N_{\bP^1/Q}\cong \cO_{\bP^1}(2)^{\oplus(r-2)}$. Once more, the short exact sequence 
\[
0\to T_{\bP^1}\to h^*T_Q\to N_{\bP^1/Q}\to 0
\]
splits, giving $h^*T_Q\cong \cO_{\bP^1}(2)^{\oplus(r-1)}$.
\end{proof}

\begin{proposition}
\label{prop:add_line}
If $\QEx(g,r,d-1)$ holds, then so does $\QEx(g,r,d)$.
\end{proposition}

\begin{proof}
Let $C'$ be a general curve of genus $g$, and let $f':C'\to Q$ be a nice map of degree $d-1$. Choose a point $x\in C'$, and let $h:\bP^1\hookrightarrow Q$ be the inclusion of a line passing through (the image of) $x$. Then, gluing along $x$ at a node yields a map $f:C'\cup\bP^1\to Q$ of degree $d$ and genus $g$. Clearly, $f$ is non-degenerate; we claim that $H^1(C'\cup\bP^1,f^*T_Q)=0$.

Indeed, the normalization map $\nu:C'\sqcup\bP^1\to C$ induces a short exact sequence
\[
0\to f^*T_Q\to \nu_*\nu^*(f^*T_Q)\to T_{Q,x}\to 0.
\]
The long exact sequence on cohomology gives
\[
\begin{aligned}
H^0(C',(f')^*T_Q)\oplus H^0(\bP^1,h^*T_Q)
&\to T_{Q,x}\to H^1(C,f^*T_Q)\\
&\to H^1(C',(f')^*T_Q)\oplus H^1(\bP^1,h^*T_Q)=0.
\end{aligned}
\]
The vanishing of the last term follows from the assumption that $f'$ is nice and Lemma \ref{lem:line_tangent}. On the other hand, the map $H^0(\bP^1,h^*T_Q)\to T_{Q,x}$ has cokernel $H^1(\bP^1,h^*T_Q(-1))=0$, also by Lemma \ref{lem:line_tangent}. Therefore, we have $H^1(C,f^*T_Q)=0$. 
\end{proof}

\begin{proposition}
\label{prop:add_conic}
If $\QEx(g-2,r,d-2)$ holds, then so does $\QEx(g,r,d)$.
\end{proposition}

\begin{proof}
Let $C'$ be a general curve of genus $g-2$, and let $f':C'\to Q$ be a nice map of degree $d-2$. Arguing as in the proof of Lemma \ref{lem:point_general}, we may choose 3 general points $x_1,x_2,x_3\in C'$, whose images lie on a smooth conic $h:\bP^1\to Q$, given by the transverse intersection of $Q$ with a 2-plane. Let $f:C=C'\cup\bP^1\to Q$ be the map of degree $d$ and genus $g$ obtained by gluing along $x_1,x_2,x_3$. Arguing using the normalization exact sequence as in Proposition \ref{prop:add_line}, and using that, by Lemma \ref{lem:conic_tangent}, $H^1(\bP^1,h^*T_Q(-3))=0$, we have $H^1(C,f^*T_Q)=0$. We conclude as in the proof of Proposition \ref{prop:add_line}.
\end{proof}

\begin{proof}[Proof of Theorem \ref{thm:main_existence}]
Fix $r\ge4$ throughout. By Proposition~\ref{prop:add_line}, for fixed $g$, it suffices to prove $\QEx(g,r,d)$ when
\begin{equation}
\label{eq:d_min}
d=\max\left\{
\left\lceil\frac{r}{r+1}\cdot g+r\right\rceil,\,
g+\left\lfloor\frac r2\right\rfloor+1
\right\}.
\end{equation}
We assume \eqref{eq:d_min} for the rest of the proof, and proceed by induction on $g$.

Write $g=q(r+1)+\epsilon$, where $q$ is an integer and $0\leq\epsilon\leq r$. For base cases, assume first that $\epsilon\in\{0,1\}$, and that $q\leq\lceil\frac r2\rceil-1$. We claim that $\binom{r+2}{2}>2d-g+1$, so we obtain $\QEx(g,r,d)$ from Proposition~\ref{prop:existence_surjective_range}.

Indeed, assume first that $\epsilon=1$, and note that
\[
\begin{aligned}
\left\lceil\frac{r}{r+1}\cdot g+r\right\rceil
&=(qr+1)+r\\
&\geq (qr+1)+\left(\left\lfloor\frac r2\right\rfloor+q+1\right)\\
&= g+\left\lfloor\frac r2\right\rfloor+1,
\end{aligned}
\]
so $d=qr+r+1$. Then, we have
\[
\begin{aligned}
2d-g+1
&=2(qr+r+1)-(qr+q+1)+1\\
&=q(r-1)+2r+2\\
&\leq \left(\frac r2-\frac12\right)(r-1)+2r+2\\
&=\frac12r^2+r+\frac52\\
&< \binom{r+2}{2}.
\end{aligned}
\]
(The last inequality requires $r\ge4$.) Replacing $\epsilon=1$ with $\epsilon=0$ decreases both $d$ and $g$ by 1, so the same inequality holds.

Now, assume that $\QEx(g,r,d)$ holds for all smaller values of $g$. If $\epsilon\geq 2$, then decreasing $g$ by 2 also decreases $d$ by 2. By Proposition~\ref{prop:add_conic}, $\QEx(g-2,r,d-2)$ implies $\QEx(g,r,d)$.

Therefore, it suffices to assume that $\epsilon\in\{0,1\}$, and that $q\geq\left\lceil\frac r2\right\rceil$. Assume that $\epsilon=0$. Then,
\[
\begin{aligned}
\left\lceil\frac{r}{r+1}\cdot g+r\right\rceil
&=qr+r\\
&\leq q(r+1)+\left\lfloor\frac r2\right\rfloor+1\\
&=g+\left\lfloor\frac r2\right\rfloor+1,
\end{aligned}
\]
and
\[
\begin{aligned}
\left\lceil\frac{r}{r+1}\cdot(g-2)+r\right\rceil
&=qr+r-1\\
&\leq q(r+1)+\left\lfloor\frac r2\right\rfloor-1\\
&=(g-2)+\left\lfloor\frac r2\right\rfloor+1.
\end{aligned}
\]
Therefore, decreasing $g$ by 2 also decreases $d$ by 2, and we conclude $\QEx(g,r,d)$ from Proposition~\ref{prop:add_conic} and the inductive hypothesis. The case $\epsilon=1$ is similar.
\end{proof}

\begin{remark}\label{rem:non_special}
Suppose that $d$ is even and that $d\ge g+r$, which narrowly falls short of the full range of Theorem \ref{thm:main_existence}. (It is typical that maps of the smallest possible degree are the most subtle to construct.) Then, one has the following direct proof of $\QEx(g,r,d)$, suggested to us by Eric Jovinelly. 

Let $f:C\to\bP^1\subset Q$ be the composition of a cover of degree $\frac{d}{2}$, followed by the inclusion of a conic spanning a 2-plane transverse to $Q$. Let $\cL$ be the pullback of $\cO_{\bP^1}(1)$ to $C$. Then, $H^1(C,\cL^{\otimes 2})=0$ by \cite[Remark 2.1]{jp2}, so
\[
H^1(C,f^*T_Q)=H^1(C,(\cL^{\otimes 2})^{\oplus(r-1)})=0,
\]
by Lemma \ref{lem:conic_tangent}.

Moreover, because
\[
h^0(C,\cL^{\otimes 2})=d-g+1\ge r+1,
\]
for general points $p_1,\ldots,p_{r+1}\in C$, the map
\[
H^0(C,f^*T_Q)\to (f^*T_Q)|_{p_1+\cdots+p_{r+1}}
\]
is surjective. It follows that $f$ deforms in a family of curves interpolating through $r+1$ general points of $Q$, so in particular deforms to a non-degenerate (and hence nice) map.

If instead $d=g+\ell$, where $2\le\ell\le r-1$, then the same construction yields maps $f:C\to Q$ with $H^1(C,f^*T_Q)=0$, but $f$ is only guaranteed to deform to a map whose linear span has dimension $\ell$. In fact, this is sharp, because $h^0(C,\cO_{\bP^r}(1))=h^0(C,\cL^{\otimes 2})=\ell+1$. This is consistent with the fact that the expected dimension of the space of maps $f:C\to Q$, factoring through the intersection $Q_\ell=Q\cap\Lambda$ of $Q$ with some $\ell$-plane, is
\[
\begin{aligned}
\dim(\Gr(\ell+1,r+1))+\chi(C,f^*T_{Q_\ell})
&=(\ell+1)(r-\ell)+(\ell-1)(d-g+1)\\
&=(r-1)(d-g+1)\\
&=\chi(C,f^*T_Q).
\end{aligned}
\]
Therefore, $f$ lies on a component of the moduli space of maps $C\to Q$ that is generically smooth of expected dimension, but parametrizes only degenerate maps. 

The same calculation shows that, if $f$ is a degenerate map with linear span of dimension greater than $d-g$, then $f$ lies on a component of the moduli space of maps $C\to Q$ of dimension strictly greater than the expected.
\end{remark}

\subsection{Strong Maximal Rank}\label{sec:smr}

Suppose that $\rho(g,r,d)\geq 0$ and $\binom{r+2}{2}\leq 2d-g+1$. By Maximal Rank, if $f:C\to\bP^r$ is a general BN-curve, then the restriction map
\[
\mu:H^0(\bP^r,\cO_{\bP^r}(2))\to H^0(C,\cO_C(2))
\]
is injective. However, if $d\geq g+\lfloor\frac{r}{2}\rfloor+1$, then maps $f:C\to Q\subset\bP^r$, proven to exist in Theorem \ref{thm:main_existence}, have the property that $\mu$ is not injective. That is, $f$ lies in the proper closed subset $Z_{g,d}^r\subset\cM_d(C,\bP^r)^\circ$ in the space of non-degenerate maps where $\mu$ fails to be injective. 

As a degeneracy locus, the expected codimension of $Z_{g,d}^r$ is
\[
(2d-g+1)-\binom{r+2}{2}+1,
\]
which is an upper bound for the actual codimension. Under additional numerical assumptions (see the discussion after Proposition \ref{prop:maps_expected_codimension} below), the \emph{Strong Maximal Rank Conjecture} of Aprodu-Farkas \cite{af} predicts that equality holds.

\begin{proposition}\label{prop:maps_expected_codimension}
With assumptions as in Theorem \ref{thm:main_existence}, suppose that $\binom{r+2}{2}\leq 2d-g+1$. Let $f:C\to Q\subset\bP^r$ be a nice map with the property that the restriction map
\[
\mu:H^0(\bP^r,\cO_{\bP^r}(2))\to H^0(C,\cO_C(2))
\]
has 1-dimensional kernel. Then, $f$ lies on a component of $Z_{g,d}^r\subset\cM_d(C,\bP^r)^\circ$ of expected codimension.
\end{proposition}

If $\binom{r+2}{2}\leq 2d-g+1$ and there exists a non-degenerate map $f:C\to Q\subset\bP^r$, then $\ker(\mu)\neq 0$. The hypothesis in Proposition \ref{prop:maps_expected_codimension} is that $Q$ is the \emph{only} quadric containing the image of $f$. In fact, in the setting of Proposition \ref{prop:maps_expected_codimension}, such nice maps always exist, as we prove in Theorem \ref{thm:existence_injective}. Because the proof proceeds by tracing back through the proof of Theorem \ref{thm:main_existence} with various technical modifications, we take Theorem \ref{thm:existence_injective} as given for the rest of this section, and relegate the proof to Appendix \S\ref{sec:quadric_unique}.

\begin{proof}[Proof of Proposition \ref{prop:maps_expected_codimension}]
As in the proof of Proposition \ref{prop:existence_surjective_range}, fix a general curve $C$, and consider the incidence correspondence
\[
\begin{tikzcd}[column sep=large]
& Y=\{(f,F)\mid \operatorname{im}(f)\subset V(F)\} \arrow[dl, "\pr_1"'] \arrow[dr, "\pr_2"] & \\
\cM_d(C,\bP^r)^\circ & & \bP H^0(\cO_{\bP^r}(2))
\end{tikzcd}
\]
of pairs of a curve contained in a quadric, where $\cM_d(C,\bP^r)^\circ$ is the space of non-degenerate maps $f:C\to\bP^r$ of degree $d$.
A nice map $f:C\to Q\subset\bP^r$, as in Theorem \ref{thm:main_existence}, defines a point of $Y$ that is smooth in the fiber $\pr_2^{-1}([Q])$, so moves in a family of dimension $\chi(C,f^*T_Q)$. Therefore, $f$ defines a smooth point of $Y$ of dimension
\[
\begin{aligned}
&\chi(C,f^*T_Q)+\binom{r+2}{2}-1\\
&=\chi(C,f^*T_{\bP^r})-\left[(2d-g+1)-\binom{r+2}{2}+1\right].
\end{aligned}
\]
If $f$ has the property that the restriction map $\mu$ has 1-dimensional kernel, then $\pr_1$ is an isomorphism near $f$. Furthermore, the fibers of $\pr_1$ are irreducible, so $\pr_1$ is an isomorphism onto its image in a neighborhood of $f\in \cM_d(C,\bP^r)^\circ$. Therefore, $f\in\cM_d(C,\bP^r)^\circ$ lies on a component of $Z_{g,d}^r=\im(\pr_1)$ of expected codimension.
\end{proof}

Suppose now that
\begin{equation}\label{eq:smr}
d\leq g+r \text{ and } 0\leq\rho(g,r,d)<r-2.
\end{equation}
Let $C$ be a fixed, general curve of genus $g$, and let $G_d^r(C)$ be the moduli space of rank $r$, degree $d$ linear series on $C$. The numerical assumptions \eqref{eq:smr} guarantee that any such linear series is very ample and complete, so may in particular be identified with the underlying line bundle $\cL\in\mathrm{Pic}^d(C)$. Let $\Sigma_{g,d}^r(C)\subset G_d^r(C)$ be the degeneracy locus where the map
\[
\nu:\operatorname{Sym}^2 H^0(C,\cL)\to H^0(C,\cL^{\otimes 2})
\]
is not of maximal rank.

Then, the Strong Maximal Rank Conjecture of Aprodu-Farkas \cite{af} predicts that $\Sigma_{g,d}^r(C)\subset G_d^r(C)$ has codimension
\[
\left|(2d-g+1)-\binom{r+2}{2}\right|+1.
\]
(It makes a similar prediction for hypersurfaces of any degree.) We refer to \cite[\S2]{robu} for a survey of results toward the Strong Maximal Rank Conjecture.

\begin{corollary}\label{cor:smr}
Let $r\ge 4$ be an integer. Suppose that $\binom{r+2}{2}\leq 2d-g+1$, that $\lfloor\frac{r}{2}\rfloor+1\leq d-g\leq r$, and that $0\leq\rho(g,r,d)<r-2$. Then, there exists a component of $\Sigma_{g,d}^r(C)\subset G_d^r(C)$ satisfying the prediction of the Strong Maximal Rank Conjecture.
\end{corollary}

\begin{proof}
Under these numerical assumptions, near a nice map $f:C\to Q\subset\bP^r$ as in Proposition \ref{prop:maps_expected_codimension}, $\cM_d(C,\bP^r)^\circ$ is an $\mathrm{Aut}(\bP^r)$-bundle over $G_d^r(C)$, and the maps $\mu,\nu$ are identified. Thus, granting Theorem \ref{thm:existence_injective}, the conclusion is immediate from Proposition \ref{prop:maps_expected_codimension}.
\end{proof}

Robu \cite[Theorems 1.1-1.3]{robu} found many examples of $\Sigma_{g,d}^r(C)\subset G_d^r(C)$ containing components of expected dimension, when either $r\leq 4$ or $d>g+r$, the latter going beyond the numerical assumptions \eqref{eq:smr} given in \cite{af}. Robu also gives examples where, when \eqref{eq:smr} does not hold, Strong Maximal Rank can fail. In a different direction, Cotterill-Alonso Gonzalo-Zhang \cite[Theorem 1.1]{cgz} established the non-emptiness of $\Sigma^r_{g,d}(C)$ in many examples, without detecting components of expected codimension. As far as we are aware, the range in which Corollary \ref{cor:smr} holds is new.

\appendix

\section{The injective range}\label{sec:quadric_unique}

The purpose of the appendix is to prove the following refinement of Theorem \ref{thm:main_existence}, which we needed to deduce Corollary \ref{cor:smr}.

\begin{theorem}
\label{thm:existence_injective}
With assumptions as in the ``if'' direction of Theorem \ref{thm:main_existence}, suppose furthermore that
\[
\binom{r+2}{2}-(2d-g+1)\leq 0.
\]
Then, there exists a nice map $f\colon C\to Q\subset\bP^r$ with the property that the restriction map
\[
\mu\colon H^0(\bP^r,\cO_{\bP^r}(2))\to H^0(C,\cO_C(2))
\]
has 1-dimensional kernel. That is, $Q$ is the unique quadric containing the image of $f\colon C\to\bP^r$.
\end{theorem}

The inequality $\binom{r+2}{2}-(2d-g+1)\leq 0$ defines the ``injective range'' in the sense that, by maximal rank, $\mu$ is injective for a general BN-curve $f:C\to\bP^r$. Theorem \ref{thm:existence_injective} shows that there exist nice maps $f:C\to Q\subset\bP^r$ failing the prediction of the Maximal Rank Conjecture in the mildest possible way.

\begin{proposition}
\label{prop:existence_injective_reduction}
It suffices to prove Theorem \ref{thm:existence_injective} under the assumption $\binom{r+2}{2}=2d-g+1$.
\end{proposition}

\begin{proof}
We will prove Theorem \ref{thm:existence_injective} assuming that it holds when $\binom{r+2}{2}=2d-g+1$. We fix $r$ and proceed by induction on $d$. Assume first that
\[
\binom{r+2}{2}-(2d-g+1)\leq -2.
\]
Then, recall from \S\ref{sec:inductive_existence} that the statement $\QEx(g,r,d)$ is deduced from either $\QEx(g,r,d-1)$ or $\QEx(g-2,r,d-2)$, by attaching to a nice map $f'\colon C'\to Q$ of lower degree either a 1-secant line or a 3-secant conic, respectively. In either case, note that the value of $\binom{r+2}{2}-(2d-g+1)$ increases by 2 upon decrementing $(g,r,d)$ to $(g,r,d-1)$ or $(g-2,r,d-2)$. Therefore, by the inductive hypothesis, we may assume that $Q$ is the unique quadric containing the image of $f'\colon C'\to\bP^r$. The same is true of the attached map $f\colon C'\cup\bP^1\to Q$, as well as a general deformation of $f$. Therefore, Theorem \ref{thm:existence_injective} holds for $(g,r,d)$.

If instead $\binom{r+2}{2}-(2d-g+1)=-1$, then the nice map $f'\colon C'\to Q$ is obtained from Proposition \ref{prop:existence_surjective_range}. In particular, $f'\colon C'\to\bP^r$ satisfies maximal rank, so $Q$ is the unique quadric containing the image of $f'$. Again, we conclude that the same is true for $f\colon C'\cup\bP^1\to Q$, so Theorem \ref{thm:existence_injective} holds for $(g,r,d)$.
\end{proof}

Therefore, to prove Theorem \ref{thm:existence_injective}, it suffices to assume that $\binom{r+2}{2}=2d-g+1$. In this case, $\QEx(g,r,d)$ is obtained from a map $f'\colon C'\to Q\subset\bP^r$ of lower degree satisfying maximal rank, so
\[
\dim (H^0(\mathcal{I}_{f'}(2)))=2.
\]
That is, the image of $f'$ lies on a 1-dimensional projective space of quadrics, including $Q$. If the attached map $f\colon C'\cup\bP^1\to Q$ has the property that $Q$ is the unique quadric containing the image of $f$, then we may conclude Theorem \ref{thm:existence_injective} for $(g,r,d)$ after taking a general deformation of $f$. Equivalently, we may conclude Theorem \ref{thm:existence_injective} if attaching the rational component $\bP^1$ to $C'$ imposes a non-trivial linear condition on quadrics containing the image of $C'$. We will prove the following more general statement, which is a refinement of Theorem \ref{thm:main_smooth}.

\begin{theorem}
\label{thm:smooth_conic}
Adopt the setup of Theorem \ref{thm:main_smooth}, and assume furthermore that
\[
\binom{r+2}{2}-(2d-g+1)\geq 2.
\]
Let $f\colon C\to\bP^r$ be a general (in particular, satisfying maximal rank) BN-curve of degree $d$ and genus $g$, and let $Q\subset\bP^r$ be a general smooth quadric containing the image of $f$.
\begin{enumerate}
\item[(a)] There exists a 1-secant line $L\subset Q$ to $C$ imposing a non-trivial condition on quadrics containing the image of $f$. That is, we have
\[
0<h^0(\mathcal{I}_{C\cup L\to \bP^r}(2))<h^0(\mathcal{I}_C(2)),
\]
where $C\cup L$ is the nodal curve obtained by gluing $C$ and $L$ along a single point of intersection, and $C\cup L\to Q$ is the resulting morphism.
\item[(b)] For general points $y_1,y_2,y_3\in C$, the unique smooth conic $K\subset Q$ containing the $y_i$ imposes a non-trivial condition on quadrics containing the image of $f$. That is, we have
\[
0<h^0(\mathcal{I}_{C\cup K}(2))<h^0(\mathcal{I}_C(2)).
\]
\end{enumerate}
\end{theorem}

Theorem \ref{thm:existence_injective} will follow, by applying either Theorem \ref{thm:smooth_conic}(a) to $(g,r,d-1)$ or Theorem \ref{thm:smooth_conic}(b) to $(g-2,r,d-2)$. (We have repurposed the letters $g,r,d$ in the statement of Theorem \ref{thm:smooth_conic} for sake of clarity of what follows.)

\begin{proof}[Proof of Theorem \ref{thm:smooth_conic}(a)]
Let $y\in Q$ be a general point, with the property that not every quadric containing the image of $C$ contains $y$. Then, the tangent hyperplane $T_yQ\subset \bP^r$ intersects $C$ non-trivially at some point $y'\in C$. Taking $L\subset Q$ to be the unique line containing $y,y'$, we have
\[
0<h^0(\mathcal{I}_{C\cup L}(2))<h^0(\mathcal{I}_C(2)).
\]
Indeed, $Q$ contains $C\cup L$, but not every quadric in $\bP H^0(\mathcal{I}_C(2))$ contains $y\in L$.
\end{proof}

The remainder of this section is devoted to the proof of Theorem \ref{thm:smooth_conic}(b). We follow the same strategy as in the proof of Theorem \ref{thm:main_smooth}, arguing by induction on $r$. We will need to weaken the numerical assumptions slightly, allowing also $g=0$ and $3\leq d\leq r-1$. In this case, $f\colon C\to\bP^r$ is a rational normal curve contained in a linear subspace of dimension $d$, and we may assume furthermore that the linear span of $C$ intersects $Q$ in a smooth quadric. We explain the inductive steps first.

\begin{proposition}
\label{prop:smooth_conic_inductive}
Assume that the conclusion of Theorem \ref{thm:smooth_conic}(b) holds when $3\leq d\leq r+2$. (If $d\leq r$, then necessarily $g=0$.) Then, Theorem \ref{thm:smooth_conic}(b) holds in general.
\end{proposition}

\begin{proof}
We follow the proof of Theorem \ref{thm:main_smooth} given in \S\ref{sec:inductive_smooth}. We proceed by induction on $r$; when $r=3$, then either $(g,d)=(0,3)$ or $(1,4)$, both of which lie in the range $3\leq d\leq r+2$. Thus, we may assume that $r\geq 4$ and that Theorem \ref{thm:smooth_conic}(b) holds for $r-1$. We may also assume that $d\geq r+3$.

We focus first on the case $d\geq g+r$ and $g\leq r-2$; the other cases are easier. Adopting the same notation as in the proof of Theorem \ref{thm:main_smooth}, fix a hyperplane $H\subset\bP^r$. Write $d''=d-r\ge 3$, let $C''$ be a smooth curve of genus 0, and let $f''\colon C''\to H$ be a general map of degree $d''$. By Theorem \ref{thm:main_smooth}, we may choose a smooth quadric $Q'\subset H$ containing the image of $f''$. If $d''\leq r-2$, then we may assume that $Q'$ intersects the linear span $\Lambda$ of $C''$ in a smooth quadric. Moreover, we may assume that the vector space of quadrics in $H$ containing the image of $f''$ has dimension
\[
h^0(\mathcal{I}_{C''}(2))=\binom{r+1}{2}-(2d''+1)=:m.
\]

By the inductive hypothesis, we may choose 3 points $y_1,y_2,y_3\in C''$ such that the unique smooth conic $K\subset Q'$ containing $y_1,y_2,y_3$ satisfies
\[
h^0(\mathcal{I}_{C''\cup K}(2))\leq m-1.
\]
Now, we choose general points $x_1,\ldots,x_{g+1}\in C''$ and $x_{g+2},\ldots,x_r\in Q'$, as in the proof of Theorem \ref{thm:main_smooth}. In particular, we may choose $x_{g+2},\ldots,x_r\in Q'$, in such a way that
\[
\begin{aligned}
h^0(\mathcal{I}_{C''\cup K\cup\{x_{g+2},\ldots,x_r\}}(2))&\leq (m-1)-(r-g-1)\\
&=\binom{r+2}{2}-(2d-g+1)-1.
\end{aligned}
\]

Finally, attaching to $C''$ a union of rational curves $h\colon R\to Q$, as in the proof of Theorem \ref{thm:main_smooth}, we obtain a map $f\colon C=C'\cup C''\to Q$.  We compute that
\[
h^0(\mathcal{I}_C(2))=h^0(\mathcal{I}_{C''\cup\{x_{g+2},\ldots,x_r\}}(2))=\binom{r+2}{2}-(2d-g+1),
\]
but
\[
h^0(\mathcal{I}_{C\cup K}(2))=h^0(\mathcal{I}_{C''\cup K\cup\{x_{g+2},\ldots,x_r\}}(2))\leq\binom{r+2}{2}-(2d-g+1)-1.
\]
Therefore, $f\colon C\to Q\subset\bP^r$ is a BN-curve for which there exist points $y_i\in C$ with the needed property asserted in Theorem \ref{thm:smooth_conic}(b), so the same holds for a general BN-curve.

In the remaining cases ($d<g+r$, or $d\geq g+r$ and $g\geq r-1$), the argument is the same, except that all of the points $x_i$ are chosen on $C''$.
\end{proof}

We now deal with the remaining cases, where $3\leq d\leq r+2$. In this range, the inductive argument of Proposition \ref{prop:smooth_conic_inductive} breaks down, because there is insufficient room to choose suitable points $y_1,y_2,y_3\in C''$. However, we will show that some of the $y_i$ can be chosen instead among the points $x_i\in H$.

\begin{lemma}
\label{lem:conic_condition_base_cases}
Suppose $r\geq 3$, let $H\subset\bP^r$ be a hyperplane, and let $Q'\subset H$ be a smooth quadric.
\begin{enumerate}
\item[(a)] Let $x_1,\ldots,x_s\in Q'$ be general points, where $s\in\{r,r+1\}$. Let $K\subset Q'$ be the unique smooth conic containing $x_1,x_2,x_3$. Then,
\[
h^0(\mathcal{I}_{K\cup\{x_4,\ldots,x_s\}\subset H}(2))<h^0(\mathcal{I}_{\{x_1,\ldots,x_s\}\subset H}(2)).
\]
That is, $K$ imposes a non-trivial condition on quadrics in $H$ passing through $x_1,\ldots,x_s$.
\item[(b)] Suppose $r\geq 4$. Let $L\subset Q'$ be a line, let $x_2,\ldots,x_{r+1}\in Q'$ be general points, and let $K\subset Q'$ be the unique smooth conic containing $x_2,x_3,x_4$. Then,
\[
h^0(\mathcal{I}_{L\cup K\cup\{x_5,\ldots,x_{r+1}\}\subset H}(2))<h^0(\mathcal{I}_{L\cup\{x_2,\ldots,x_{r+1}\}\subset H}(2)).
\]
\item[(c)] Suppose $r\geq 4$. Let $K_0\subset Q'$ be a smooth conic, whose linear span is not contained in $Q'$. Let $g\in\{0,1,2\}$ be an integer, let $x_{g+2},\ldots,x_{r}\in Q'$ and $y_1,y_2\in K_0$ be general points, and let $K\subset Q'$ be the unique smooth conic containing $y_1,y_2,x_{g+2}$. Then,
\[
h^0(\mathcal{I}_{K_0\cup K\cup\{x_{g+3},\ldots,x_{r}\}\subset H}(2))<h^0(\mathcal{I}_{K_0\cup\{x_{g+2},\ldots,x_{r}\}\subset H}(2)).
\]
\end{enumerate}
\end{lemma}

The unusual indexing in parts (b) and (c) is meant to match the notation in the proof of Theorem \ref{thm:main_smooth}, in the case $d\geq g+r$ and $g\leq r-2$.

\begin{proof}
We begin with (a). We have
\[
h^0(\mathcal{I}_{K\subset H}(2))=\binom{r+1}{2}-5<\binom{r+1}{2}-3=h^0(\mathcal{I}_{\{x_1,x_2,x_3\}\subset H}(2)).
\]
Furthermore, the remaining points $x_4,\ldots,x_s$ are general, so impose independent linear conditions on each linear system, until possibly shrinking them to the 1-dimensional vector space of forms cutting out $Q'$. On the other hand,
\[
h^0(\mathcal{I}_{\{x_1,\ldots,x_s\}\subset H}(2))=\binom{r+1}{2}-s>1,
\]
so imposing the additional conditions at $x_4,\ldots,x_s$ does not shrink both linear systems to the 1-dimensional vector space of forms cutting out $Q'$. The conclusion follows.

Next, we prove (b). We claim first that
\[
h^0(\mathcal{I}_{L\cup K\subset H}(2))<h^0(\mathcal{I}_{L\cup\{x_2,x_3,x_4\}\subset H}(2)).
\]
Indeed, because $x_2,x_3,x_4\in Q'$ are general, the right-hand side is equal to $\binom{r+1}{2}-6$. If $r=4$, then $L\cup K\subset Q'\cong\bP^1\times\bP^1$ is a divisor of bidegree $(1,2)$. Thus, $H^0(\mathcal{I}_{L\cup K\subset H}(2))\subset H^0(\bP^3,\cO(2))$ is identified with the pullback of
\[
H^0(\bP^1\times\bP^1,\cO(1,0))\subset H^0(\bP^1\times\bP^1,\cO(2,2))
\]
under the surjection $H^0(\bP^3,\cO(2))\to H^0(\bP^1\times\bP^1,\cO(2,2))$ with 1-dimensional kernel. In particular, we have
\[
h^0(\mathcal{I}_{L\cup K\subset H}(2))=3<\binom{4+1}{2}-6,
\]
as needed.

If $r\geq 5$, then because $x_2,x_3,x_4\in Q'$ are general, $L,K$ impose independent linear conditions on quadrics in $H$. Indeed, we may choose coordinates $z_1,\ldots,z_r$ on $H$ in such a way that $L$ imposes 3 independent conditions on the monomials in $z_1,z_2$, and $K$ imposes 5 independent conditions on the monomials in $z_3,z_4,z_5$. In particular, we have
\[
h^0(\mathcal{I}_{L\cup K\subset H}(2))=\binom{r+1}{2}-8.
\]
Now, the remaining points $x_5,\ldots,x_{r+1}$ impose independent linear conditions on both linear systems. Furthermore, we have
\[
h^0(\mathcal{I}_{L\cup\{x_2,\ldots,x_{r+1}\}\subset H}(2))=\binom{r+1}{2}-(r+3)>0,
\]
so (b) follows.

Finally, consider (c). We claim first that
\[
h^0(\mathcal{I}_{K_0\cup K\subset H}(2))<h^0(\mathcal{I}_{K_0\cup\{x_{g+2}\}\subset H}(2)).
\]
The needed statement will follow from a similar argument as in the previous two cases. The right-hand side is equal to $\binom{r+1}{2}-6$; we claim that the union of conics $K_0\cup K$ imposes 8 conditions on quadrics in $H$. It suffices to verify this claim for quadrics inside the linear span $\Lambda\cong\bP^3$ of $K_0\cup K$, by the surjectivity of the restriction map $H^0(H,\cO(2))\to H^0(\Lambda,\cO(2))$. 

By the generality assumptions, we may assume that $\Lambda\cap Q'\cong\bP^1\times\bP^1$ is smooth. Then, $H^0(\mathcal{I}_{K_0\cup K\subset\Lambda}(2))$ is identified with the pullback of
\[
H^0(\bP^1\times\bP^1,\cO)\subset H^0(\bP^1\times\bP^1,\cO(2,2))
\]
under $H^0(\bP^3,\cO(2))\to H^0(\bP^1\times\bP^1,\cO(2,2))$. In particular, we have $h^0(\mathcal{I}_{K_0\cup K\subset\Lambda}(2))=2$, as needed.
\end{proof}

\begin{corollary}
\label{cor:smooth_conic_mid_degree}
Theorem \ref{thm:smooth_conic}(b) holds in the following cases.
\begin{enumerate}
\item[(a)] $(g,d)=(0,r)$ or $(1,r+1)$.
\item[(b)] $(g,d)=(0,r+1)$.
\item[(c)] $d=r+2$ and $g\leq 2$.
\end{enumerate}
\end{corollary}

\begin{proof}
We apply the same arguments as in the proof of Proposition \ref{prop:smooth_conic_inductive}, where the conic $K$ is now taken to be the conic from the corresponding part of Lemma \ref{lem:conic_condition_base_cases}. In the case $(g,d)=(0,r)$, let $x_1,\ldots,x_r\in Q'$ be general points as in Lemma \ref{lem:conic_condition_base_cases}(a), and let $h:R\to Q$ be as in Proposition \ref{prop:existence_chain}. Then, $h$ is a BN-curve, and taking $y_i=x_i$ for $i=1,2,3$ gives the desired conic. 

In the case $(g,d)=(1,r+1)$, we take $x_1,\ldots,x_{r+1}\in Q'$, and (as in the proof of Theorem \ref{thm:main_smooth}) $h:E\to Q'$ to be as in Proposition \ref{prop:existence_elliptic}. We proved Proposition \ref{prop:existence_elliptic} only when $r\geq 4$, but when $r=3$, we have the following direct argument. A general BN-curve $f:E\to\bP^3$ of degree 4 has $h^0(\cI_E(2))=2$. The base locus of the corresponding pencil of quadrics has degree 4, so must be $E$ itself. In particular, any 3-secant conic to $E$ imposes a nontrivial condition on this linear system.

If $(g,d)=(0,r+1)$, then adopting the notation of the proof of Theorem \ref{thm:main_smooth}, we take $f'':C''\to Q'$ to be the inclusion of the line $L$ from Lemma \ref{lem:conic_condition_base_cases}(b), and again $h:R\to Q$ as in Proposition \ref{prop:existence_chain}. Then, Lemma \ref{lem:conic_condition_base_cases}(b) shows that the conic $K\subset Q'$ again has the needed property. Finally, the case $d=r+2$ is similar, where we take $f''$ instead to be the inclusion of the conic $K_0$ from Lemma \ref{lem:conic_condition_base_cases}(c).
\end{proof}

\begin{proof}[Proof of Theorem \ref{thm:smooth_conic}(b)]
Combining Proposition \ref{prop:smooth_conic_inductive} and Corollary \ref{cor:smooth_conic_mid_degree}, it remains to consider the case $g=0$ and $3\leq d<r$. In this case, Theorem \ref{thm:smooth_conic}(b) follows from the case $(g,r,d)=(0,d,d)$. Indeed, there exists a rational normal curve $f\colon C\hookrightarrow\bP^d$, a smooth quadric $Q'\subset\bP^d$ containing $C$, and points $y_1,y_2,y_3\in C$ such that the unique conic $K\subset Q'$ containing the $y_i$ imposes a non-trivial condition on $H^0(\mathcal{I}_{C\subset\bP^d}(2))$. Now, let $\bP^d\subset\bP^r$ be the inclusion of a linear subspace, and choose a lift $Q\subset\bP^r$ of $Q'$ such that $Q\cap\bP^d=Q'$. Then, $K\subset Q$ is still the unique conic containing $y_1,y_2,y_3$, and imposes a non-trivial condition on $H^0(\mathcal{I}_{C\subset\bP^r}(2))$.
\end{proof}

\bibliographystyle{alpha} 
\bibliography{existence_quadric_v2.bib}

\newcommand{\etalchar}[1]{$^{#1}$}
\begin{thebibliography}{LOTiBZ24}

\bibitem[ACGH85]{acgh}
E.~Arbarello, M.~Cornalba, P.~A. Griffiths, and J.~Harris.
\newblock {\em Geometry of algebraic curves. {V}ol. {I}}, volume 267 of {\em
  Grundlehren der mathematischen Wissenschaften [Fundamental Principles of
  Mathematical Sciences]}.
\newblock Springer-Verlag, New York, 1985.

\bibitem[AF11]{af}
M.~Aprodu and G.~Farkas.
\newblock Koszul cohomology and applications to moduli.
\newblock In {\em Grassmannians, moduli spaces and vector bundles}, volume~14
  of {\em Clay Math. Proc.}, pages 25--50. Amer. Math. Soc., Providence, RI,
  2011.

\bibitem[ALY19]{aly}
A.~Atanasov, E.~Larson, and D.~Yang.
\newblock Interpolation for normal bundles of general curves.
\newblock {\em Mem. Amer. Math. Soc.}, 257(1234):v+105, 2019.

\bibitem[Bal12]{ballico}
E.~Ballico.
\newblock Embeddings of general curves in projective spaces: the range of the
  quadrics.
\newblock {\em Lith. Math. J.}, 52(2):134--137, 2012.

\bibitem[BF18]{bafa}
B.~Bakker and G.~Farkas.
\newblock The {M}ercat conjecture for stable rank 2 vector bundles on generic
  curves.
\newblock {\em Amer. J. Math.}, 140(5):1277--1295, 2018.

\bibitem[BLL{\etalchar{+}}24]{bllrst}
R.~Beheshti, B.~Lehmann, C.~Lian, E.~Riedl, J.~Starr, and S.~Tanimoto.
\newblock On the asymptotic enumerativity property for {F}ano manifolds.
\newblock {\em Forum Math. Sigma}, 12:Paper No. e112, 2024.

\bibitem[BP21]{bp}
A.~Buch and R.~Pandharipande.
\newblock Tevelev degrees in {G}romov-{W}itten theory.
\newblock {\em arXiv 2112.14824}, 2021.

\bibitem[CAGZ21]{cgz}
E.~Cotterill, A.~Alonso~Gonzalo, and N.~Zhang.
\newblock The strong maximal rank conjecture and higher rank {B}rill-{N}oether
  theory.
\newblock {\em J. Lond. Math. Soc. (2)}, 104(1):169--205, 2021.

\bibitem[CFKM14]{cfkm}
I.~Ciocan-Fontanine, B.~Kim, and D.~Maulik.
\newblock Stable quasimaps to {GIT} quotients.
\newblock {\em J. Geom. Phys.}, 75:17--47, 2014.

\bibitem[CL26]{cl_secant}
A.~Cela and C.~Lian.
\newblock Brill-{N}oether existence with secants, 2026.
\newblock arXiv:2609.26713.

\bibitem[EH83]{eh_petri}
D.~Eisenbud and J.~Harris.
\newblock A simpler proof of the {G}ieseker-{P}etri theorem on special
  divisors.
\newblock {\em Invent. Math.}, 74(2):269--280, 1983.

\bibitem[EH87a]{eh_irreducible}
D.~Eisenbud and J.~Harris.
\newblock Irreducibility and monodromy of some families of linear series.
\newblock {\em Ann. Sci. \'Ecole Norm. Sup. (4)}, 20(1):65--87, 1987.

\bibitem[EH87b]{eh_kodaira}
D.~Eisenbud and J.~Harris.
\newblock The {K}odaira dimension of the moduli space of curves of genus $\ge
  23$.
\newblock {\em Invent. Math.}, 90(2):359--387, 1987.

\bibitem[FJP25]{fjp2}
G.~Farkas, D.~Jensen, and S.~Payne.
\newblock The {K}odaira dimension of {$\overline{\mathcal{M}}_{22}$} and
  {$\overline{\mathcal{M}}_{23}$}.
\newblock {\em Camb. J. Math.}, 13(3):431--607, 2025.

\bibitem[FR20]{fr}
G.~Farkas and R.~Rim\'anyi.
\newblock Quadric rank loci on moduli of curves and {$K3$} surfaces.
\newblock {\em Ann. Sci. \'Ec. Norm. Sup\'er. (4)}, 53(4):945--992, 2020.

\bibitem[GH80]{gh}
P.~Griffiths and J.~Harris.
\newblock On the variety of special linear systems on a general algebraic
  curve.
\newblock {\em Duke Math. J.}, 47(1):233--272, 1980.

\bibitem[Gie82]{gieseker}
D.~Gieseker.
\newblock Stable curves and special divisors: {P}etri's conjecture.
\newblock {\em Invent. Math.}, 66(2):251--275, 1982.

\bibitem[Har84]{harris_odd}
J.~Harris.
\newblock On the {K}odaira dimension of the moduli space of curves, {II}. {T}he
  even genus case.
\newblock {\em Invent. Math.}, 75(3):437--466, 1984.

\bibitem[HM82]{hm}
J.~Harris and D.~Mumford.
\newblock On the {K}odaira dimension of the moduli space of curves.
\newblock {\em Invent. Math.}, 67(1):23--88, 1982.

\bibitem[JP16]{jp2}
D.~Jensen and S.~Payne.
\newblock Tropical independence {II}: {T}he maximal rank conjecture for
  quadrics.
\newblock {\em Algebra Number Theory}, 10(8):1601--1640, 2016.

\bibitem[Kad19]{kadikoylu}
\.I. Kadik\"oyl\"u.
\newblock Variety of singular quadrics containing a projective curve.
\newblock {\em Ann. Inst. Fourier (Grenoble)}, 69(4):1879--1896, 2019.

\bibitem[KP01]{kp}
B.~Kim and R.~Pandharipande.
\newblock The connectedness of the moduli space of maps to homogeneous spaces.
\newblock In {\em Symplectic geometry and mirror symmetry ({S}eoul, 2000)},
  pages 187--201. World Sci. Publ., River Edge, NJ, 2001.

\bibitem[Lar16]{larson}
E.~Larson.
\newblock Interpolation for restricted tangent bundles of general curves.
\newblock {\em Algebra Number Theory}, 10(4):931--938, 2016.

\bibitem[Lar17]{larson_mr}
E.~Larson.
\newblock The maximal rank conjecture.
\newblock 2017.
\newblock arXiv:1711.04906.

\bibitem[Lar22]{larson_reducible}
E.~Larson.
\newblock Constructing reducible {B}rill-{N}oether curves.
\newblock {\em Doc. Math.}, 27:1953--1984, 2022.

\bibitem[Lia25]{lian_hyp}
C.~Lian.
\newblock Asymptotic geometric {T}evelev degrees of hypersurfaces.
\newblock {\em Michigan Math. J.}, 75(1):119--139, 2025.

\bibitem[LOTiBZ21]{lotz}
F.~Liu, B.~Osserman, M.~Teixidor~i Bigas, and N.~Zhang.
\newblock Limit linear series and ranks of multiplication maps.
\newblock {\em Trans. Amer. Math. Soc.}, 374(1):367--405, 2021.

\bibitem[LOTiBZ24]{lotz2}
F.~Liu, B.~Osserman, M.~Teixidor~i Bigas, and N.~Zhang.
\newblock The strong maximal rank conjecture and moduli spaces of curves.
\newblock {\em Algebra Number Theory}, 18(8):1403--1464, 2024.

\bibitem[LP25]{lian_pand}
C.~Lian and R.~Pandharipande.
\newblock Enumerativity of virtual {T}evelev degrees.
\newblock {\em Ann. Sc. Norm. Super. Pisa Cl. Sci.}, 26(1):71--89, 2025.

\bibitem[LV23]{lv}
E.~Larson and I.~Vogt.
\newblock Interpolation for {B}rill-{N}oether curves.
\newblock {\em Forum Math. Pi}, 11:Paper No. e25, 2023.

\bibitem[Mio25]{mioranci}
L.~Mioranci.
\newblock Restricted tangent bundle of rational curves on projective
  hypersurfaces, 2025.
\newblock arXiv:2507.13927.

\bibitem[MS25]{ms}
A.~Marian and S.~Sinha.
\newblock A short way of counting maps to hypersurfaces in {G}rassmannians,
  2025.
\newblock arXiv:2506.10593.

\bibitem[Rob26]{robu}
V.~Robu.
\newblock The stratified and the strong maximal rank conjecture in
  {$\mathbb{P}^3$} and {$\mathbb{P}^4$}, 2026.
\newblock arXiv:2605.24612.

\bibitem[Sev15]{severi}
F.~Severi.
\newblock Sulla classificazione delle curve algebriche e sul teorema di
  esistenza di {Riemann}.
\newblock {\em Rend. R. Acc. Naz. Lincei}, 24(5):887--888, 1915.

\bibitem[TiB03]{teixidor}
M.~Teixidor~i {B}igas.
\newblock Injectivity of the symmetric map for line bundles.
\newblock {\em Manuscripta Math.}, 112(4):511--517, 2003.

\end{thebibliography}

\end{document}